\documentclass[12pt,a4paper]{amsart}
\usepackage[T1]{fontenc}   % needed for the Polish ogonek: Gaw\k{e}dzki

\makeatletter
\newcommand{\linethrough}{\mathpalette\@thickbar}
\newcommand{\@thickbar}[2]{{#1\mkern0mu\vbox{
    \sbox\z@{$#1#2\mkern-1.5mu$}%
    \dimen@=\dimexpr\ht\tw@-\ht\z@+2\p@\relax % The +2 represents the vertical shift of the line.
    \hrule\@height0.5\p@ % The 0.5 represent the thickness on the line.
    \vskip\dimen@
    \box\z@}}
}
\makeatother
\newcommand{\mathstrike}[1]{\ensuremath{\linethrough{#1}}}

\usepackage{amsmath,amsfonts,amssymb,bbm}
\usepackage{DLdef1}

\makeatletter
\@ifundefined{th@thm}{\newtheoremstyle{thm}{2mm}{2mm}{\itshape}{}{\bf}{.}{ }{}}{}
\theoremstyle{thm}
\@ifundefined{Theorem}{}{}
\@ifundefined{Lemma}{}{}
\@ifundefined{Statement}{}{}
\@ifundefined{Proposition}{}{}
\@ifundefined{Corollary}{}{}
\theoremstyle{definition}
\@ifundefined{Remark}{\newtheorem*{Remark}{\secno Remark}}{}
\@ifundefined{Example}{}{}
\@ifundefined{rank}{\DeclareMathOperator{\rank}{rank}}{}
\@ifundefined{IM}{\DeclareMathOperator{\IM}{Im}}{}
\@ifundefined{Llra}{\newcommand{\Llra}{\Longleftrightarrow}}{}
\makeatother
\usepackage{cancel}

\newcommand{\Mod} {{\sf Mod}}
\newcommand{\Size} {\text{Size}}
\makeatletter
\@ifundefined{ScommSalgs}{}{\let\ScommSalgs\@undefined}
\makeatother
\newcommand{\ScommSalgs}{{\sf SCommSAlgs}}

\newcommand{\bean}{\begin{eqnarray*}}
\newcommand{\eean}{\end{eqnarray*}}

\newcommand{\bea}{\begin{eqnarray}}
\newcommand{\eea}{\end{eqnarray}}

\usepackage{hyperref}

 \def\mmat #1,#2,#3,#4,{\text{\small\arraycolsep=3pt $
\begin{pmatrix}#1&#2\\#3&#4\end{pmatrix}$}}

\newcommand{\del}{\partial}

\let\ssec\subsection

\renewcommand {\ssbegin}[2][*]
 {\refstepcounter{subsection}%
\if#1*
\addcontentsline{toc}{subsection}{\thesubsection.\hskip 1pc #2}%
\else
\addcontentsline{toc}{subsection}{\thesubsection.\hskip 1pc #2. #1}%
\fi
 \def \secno {\gdef \secno {}{\ssecfont
\thesubsection.\hskip 2ex}%
 }%
 \begin{#2}}

\renewcommand {\sssbegin}[2][*]
 {\refstepcounter{subsubsection}%\label{sss#1}
\if#1*
\addcontentsline{toc}{subsubsection}{\thesubsubsection.\hskip 1pc #2}%
\else
\addcontentsline{toc}{subsubsection}{\thesubsubsection.\hskip 1pc #2. #1}
\fi
 \def \secno {\gdef \secno {}{\ssecfont \thesubsubsection.\hskip 2ex}%
 }%
 \begin{#2}}

\renewcommand {\parbegin}[2][*]
 {\refstepcounter{paragraph}%\label{ssss#1}
\if#1*
\addcontentsline{toc}{paragraph}{\theparagraph.\hskip 1pc #2}%
\else
\addcontentsline{toc}{paragraph}{\theparagraph.\hskip 1pc #2. #1}
\fi
 \def \secno {\gdef \secno {}{\ssecfont \theparagraph.\hskip 2ex}%
 }%
 \begin{#2}}

\title[On odd parameters in geometry]{On odd parameters in geometry}
\author[Leites]{Dimitry Leites}

\subjclass[2020]{Primary 58A50, 17B60}

\address{Department of Mathematics, Stockholm University,
Albanov\"agen 28, SE-114 19 Stockholm, Sweden; dimleites@gmail.com}

\keywords {Simple Lie superalgebra, deformation, non-split supermanifold}
\begin{document}

% amsart requires the abstract BEFORE \maketitle
\begin{abstract} 1) In 1976, looking at the list of simple finite-dimensional complex Lie superalgebras, J.~Bernstein and I, and independently M.~Duflo, observed that some divergence-free vectorial Lie superalgebras  have deformations with odd parameters and conjectured that no other simple Lie superalgebras have such deformations.  
Here, I prove this conjecture and overview the known classification of simple finite-dimensional complex Lie superalgebras, their presentations, realizations, and relations with simple Lie (super)algebras over fields of positive characteristic.

2) Any ringed space  of the form (a manifold $M$, the sheaf of sections of the exterior algebra of a vector bundle over $M$) is called split supermanifold.  Gaw\k{e}dzki (1977) and  Batchelor (1979) proved that every smooth supermanifold is split. In 1982--83, P.~Green and Palamodov showed that a~complex-analytic supermanifold can be non-split. So far, researchers considered only even obstructions to splitness, and hence concluded that any supermanifold of superdimension $m|1$ is split. I'll show  there are non-split supermanifolds of superdimension $m|1$ \textbf{over a~base with odd parameters}, although every $m|1$-dimensional supermanifold over $\text{Spec}\ \mathbb{C}$ is split; e.g., certain  superstrings, the obstructions to their splitness depend on odd parameters.
\end{abstract}

\maketitle

%\date{Received May 1, 2006}

\section{Introduction}

In this paper, I mainly consider two types of occurrences of odd parameters:

(i) in deformations of several simple Lie superalgebras;

(ii) in obstructions to splitness of complex-analytic supermanifolds. 

The Lie superalgebras considered are mainly finite-dimensional over $\Cee$, unless otherwise specified; the characteristic $p$ of the ground field $\Kee$ can not be confused with parity also denoted by~$p$. 

I overview many related results and offer several intriguing open problems.

\ssec{Prologue: Four types of main characters, two questions} In the late 1960s, many mathematicians (first, in Moscow, and later in the whole Universe, until the vogue changed) were discussing two topics: 

(1) \textbf{Kac's classification of simple $\Zee$-graded Lie algebras of polynomial growth over~$\Cee$}  (under a~technical assumption later dismissed by O.~Mathieu in a~difficult paper \cite{M}) consisting of finite-dimensional examples (classified by Cartan and Killing), vectorial Lie algebras with polynomial coefficients (classified by Cartan) and two new ingredients: loop algebras with values in simple finite-dimensional algebras and twisted versions of these loops, and the complexified Lie algebra of vector fields on the circle, see~ \cite{K}.

 (2) \textbf{The Kostrikin-Shafarevich conjecture} on classification of simple finite-dimensional \textit{modular} (i.e., over $\Kee$ of characteristic $p>0$) Lie algebras over algebraically closed fields of characteristic $p>7$ (initially formulated only for restricted algebras); for a~review and update, see \cite{BGLLS}.  
 
 In both classifications, the $\Zee$-graded infinite-dimensional simple Lie algebras of vector fields over $\Cee$ play a~prominent role. 
 
 It became more and more obvious with time that two digressions from the main goal (classification of \textbf{simple and $\Zee$-graded} objects) are important: classification of (a) deformations of $\Zee$-graded Lie algebras and (b) Lie algebras of derivations and central extensions of simple Lie algebras.

The deformation of the Lie algebra of Hamiltonian vector fields induced by quantization of the Poisson Lie algebra leads out of the set of $\Zee$-graded Lie algebras considered in topic~ (1), and was ignored by mathematicians interested in classifications (1) and (2), whereas all other examples of simple $\Zee$-graded Lie algebras are rigid. (This claim on rigidity is correct from a~\textit{certain point of view} only; I will not digress into this mysterious area here; let me just mention an \textit{ostensible} contradiction between,  on one hand, the correct, under the definitions considered, proof of the fact that $H^2(\fg;\fg)=0$ for the loop algebras $\fg$, see \cite{LR}, and, on the other hand, the explicit examples of such deformations in \cite{GoHo1, Ho2} later evolved into a~rich theory of Krichever-Novikov  algebras, see a~review \cite{Shei}. There is no contradiction, actually, because   $H^2(\fg;\fg)=0$ does not guarantee \textbf{rigidity}, which is lack of deformations. Observe also that $H^2(\fg;\fg)\neq 0$ does not imply that the deformed algebra is not isomorphic to the initial algebra; this phenomenon is called \textit{semi-trivial deformation}, see \cite{BGL1, BLLS} with many examples.) 

Although in topic (2), i.e., in the non-restricted KSh-conjecture, the deformations come to the foreground, they still did not draw as much attention as I think they should have, except for the filtered deformations. (A plausible, very probable, reason for this negligence: of the two types of simple Lie algebras, the non-filtered \textit{deforms} --- results of deformations --- of the vectorial Lie algebras are isomorphic to the known simple Lie algebras, whereas the Lie algebras of the other type --- the ones with Cartan matrix and their simple relatives --- are rigid, at least for~$p>3$. The situation is totally different if $p=2$, see Shchepochkina's exceptional example --- the deform of the divergence-free algebra --- in \cite{BGLLS}.)

Observe two more types of characters that came to the foreground in the mid-1970s:

(3) \textbf{Lie superalgebras} were known in topology and deformation theory under the self-contradictory\footnote{The old name is self-contradictory because Lie superalgebras and Lie algebras (graded or not) satisfy different identities.}  name ``graded Lie algebras'' since at least 1940, mainly over finite fields or $\Zee$  (in which case they are super Lie rings). From the mid-1970s, Lie superalgebras started to draw attention of researchers: first, in characteristic 0, thanks to applications of Lie superalgebras to high energy physics and solid body physics. 

Some  physicists, more clairvoyant than other, began to perceive the usefulness --- in physics --- of what is nowadays called ``Lie superalgebras'' even before the notion ``Lie superalgebra'' was explicitly distinguished and separated from  ``$\Zee/2$-graded Lie algebra'', see prefaces in \cite{DSB}. Encouraged by one such physicist, Kac gave examples of three series of simple vectorial Lie superalgebras, see \cite{K0} and more accessible Kac's preface in \cite{DSB}.

(4)   For decades the mathematicians interested in simple Lie (super)algebras mostly ignored the examples of \textbf{filtered simple Lie algebras which can not be obtained as a~result of deformation of  a~$\Zee$-graded simple Lie algebra}. In physics, such Lie (super)algebras appeared in higher spin theories, see \cite{Vas, Vas1}; in mathemtics, they were used to interpret known orthogonal polynomials, see \cite{LSer, Ser}. Among such simple Lie (super)algebras, the simplest to describe is $\fsl(\lambda)$, where $\lambda\in \Cee\Pee^1$, the Lie algebra of traceless ``matrices of complex size $\lambda\in \Cee$ and even $\lambda\in \Cee\Pee^1$''.  Recall its construction.

Let $\mathfrak{g}=\mathfrak{sl}(2)$,  and $c_1=\lambda^2-1$, where $\lambda=\mu+1\in\mathbb{C}$ for the highest weight $\mu$ of an irreducible $\mathfrak{sl}(2)$-module. Let $C_1$ be the quadratic Casimir operator. The associative algebra $U_\lambda:=U(\mathfrak{g})/(C_1-c_1)$, where $U(\mathfrak{g})$ is the universal enveloping algebra of $\mathfrak{g}$ and $(C_1-c_1)$ is a two-sided ideal generated by $C_1-c_1$, has an ideal
\[I_\lambda=\begin{cases}0&\text{if $\lambda\not\in \mathbb{Z}\setminus\{0\}$}\\
\text{of finite codimension}&\text{otherwise}.\\
\end{cases}
\]

Let $A$ be any associative algebra, let $A^L$ be the Lie algebra whose space is $A$ but multiplication is given by the commutator $[a,b]:=ab-ba$ for any~$a,b\in A$.  If $\lambda=n\in \mathbb{Z}\setminus\{0\}$, then $U_n/I_n\simeq \text{Mat}(|n|)$, hence the name $\fgl(\lambda):=(U_\lambda/I_\lambda)^L$ for any $\lambda\in\mathbb{C}$, see \cite{Fei}.  

While the trace can be normalized at will,  J.~Bernstein suggested a~ definition of the trace on $\fgl(\lambda)$ such that $\tr(1)=\lambda$ for $\lambda\in(\mathbb{C}\setminus\mathbb{Z})$, see \cite{KhM}, resembling $\tr(1_{|n|})=|n|$ for $n\in\Zee\setminus \{0\}$.  Let
\begin{equation}\label{sllambda}
\fsl(\lambda):=\{X\in \fgl(\lambda)\mid \tr X=0\}.
\end{equation}

If $\lambda\not\in \mathbb{Z}\setminus\{0\}$, then $U_\lambda$ and its limit as $\lambda\to\infty\in \Cee\Pee^1$, see \cite{GL1}, are central simple as well as $U_n/I_n$ for $n\in \mathbb{Z}\setminus\{0\}$. 
According to a~Herstein--Montgomery theorem, for any central simple associative algebra $A$ with center $Z$, the Lie algebra $L(A):=(A^L)'/((A^L)'\cap Z)$, where $\fg'$ designates the first derived (commutant) of the Lie (super)algebra $\fg$, is  simple, and the same is true for Lie superalgebras under a~certain condition, see \cite[Th.3.8]{Mon}.

The Lie algebras $(U_\lambda)^L$ are ``quantized'' versions of the Lie algebras of functions on the orbits of the co-adjoint representation of simple finite-dimensional Lie algebras over $\mathbb{C}$, see \cite{Kon}, considered with the Poisson bracket. The natural filtration of $U(\fsl(2))$ induces a~filtration in $(U_\lambda)^L$; the associated graded Lie algebra is isomorphic for $\lambda\neq 0$ to the Lie algebra $\fpo(2)_\ev$ realized on the space of even degree polynomials in 2 indeterminates with respect to the Poisson bracket, see \cite{HKP}; the elements of non-negative degree in the Lie grading of $\fpo(2)_\ev$ form a~Lie algebra with codimension-3 radical complementing~$\fsl(2)$.

Consider the following generalizations of $\fsl(\lambda)$, see eq.~\eqref{sllambda}. Let $\mathfrak{g}$ be any simple finite-dimensional Lie algebra over $\mathbb{C}$, let $\rho$ be the half-sum of positive roots of $\mathfrak{g}$, let $\mu$ be the highest weight of an irreducible $\mathfrak{g}$-module with highest weight vector $v^\mu$, let the $C_i$ be the Casimirs, i.e., the generators of the center of $U(\mathfrak{g})$. Set (the shift by $\rho$ is need to get the $c_i$ symmetric relative the Weyl group)
\begin{equation}\label{U_c}
U_c:=U(\mathfrak{g})/(C_i-c_i)_{i=1}^{\text{rk} \mathfrak{g}}, \text{~~where $c_i=C_i|_{v^{\mu+\rho}}\in\Cee$}.
\end{equation}  
For almost all values of $c=(c_1,\dots, c_{\text{rk}\mathfrak{g}})$, the algebra $U_c$ is central simple, so $L(U_c)$ is a~simple Lie algebra, see \cite{Lq}, where Lie superalgebras $\fg$ in this construction are also considered.

\textbf{Question 1}: What are the deformations of the simple $\Zee$-graded  Lie superalgebras (finite-dimensional  and of polynomial growth) that can not be captured by classification of equivalence classes  of the tensors the algebra preserves, cf. \cite[Subsection~24.2]{BGLLS}?  In \S~\ref{defProof}, I will show that no simple \textbf{finite-dimensional} vectorial simple Lie superalgebra over $\Cee$ has such deformations. 

Contrariwise, the deforms of $\fosp(4|2)$ (see Subsection~\ref{sssBilMatr}) provide us with examples of Lie superalgebras that can not be described by equivalence classes of the bilinear forms preserved by certain members of the parametric family of these algebras.

\textbf{Question 2}: Over $\Cee$, is there a~ simple finite-dimensional Lie superalgebra whose associated $\Zee$-graded algebra has no negative part, like $\fsl(\lambda)$? For the answer, see \cite{K2} and the gist of the ideas of the proof in \S~\ref{Kac}: no such Lie superalgebras.  The answer for the simple finite-dimensional Lie algebras over fields of characteristic~${>3}$ is the same. The same question for Lie superalgebras over fields of positive characteristic and Lie algebras in characteristics 3 and 2 is an \textbf{Open problem}.

\ssec{Disclaimer}\label{ssDiscl} It was clear from the very beginning that certain ``relatives'' of some of the  simple Lie (super)algebras~$\fg$, such as  Lie (super)algebras of derivations and central extensions
of~$\fg$, are more (at least, no less) interesting in applications than $\fg$ itself: e.g., affine Kac--Moody (super)algebras vs. loop (super)algebras and their twisted versions, Poisson Lie (super)algebras vs. Hamiltonian Lie (super)algebras. (Observe that the affine Kac--Moody algebras are most known examples of  recently distinguished notion of \textit{double extensions}, this notion simultaneously requires a~derivation, a~central extension, and a~ third ingredient --- a~non-degenerate $\fg$-invariant symmetric bilinear form on $\fg$, see \cite{BLS} and references therein.)

Until now, the vectorial Lie superalgebras (same as vectorial Lie algebras over fields of positive characteristic) and several other types of Lie (super)algebras are often described as if the author of the description is obsessed with \textbf{simple} Lie (super)algebras. This approach leads to the notation with a~wrong emphasis that manifestly has to be improved: the \textit{natural} and \textit{basic} object is the \textit{Cartan prolong} (the result of the Cartan prolongation) rediscovered by Tanaka and generalized by Shchepochkina (see \cite{Shch,BGLLS}), whereas the simple object is the derived algebra or a~quotient of the prolong. Similar arguments identify $\fg(A)$, $\fpe(n)$, $\fq(n)$, $\fpo(0|n)$, defined in what follows, as the protagonists, more natural than their simple subquotients. 

However, the \textbf{simple} Lie (super)algebras are the first to be classified, of course.

\ssec{Simple Lie superalgebras. First definitions and examples}\label{clumsy} First examples of Lie superalgebras, actually Lie super \textbf{rings} over $\Zee$, appeared in 1941, in topology as the sets of homotopy groups with the Whitehead product, see \cite{Wh}. Associated with these examples were modular Lie superalgebras over finite fields. Similarity of Lie superalgebras and modular Lie algebras is so striking that sometimes one hears and reads that \lq
\lq when $p=2$, there is no difference between Lie algebras and Lie superalgebras'', which is not true as is lucidly explained in \cite{BGLLS}.

Weibel writes: ``Homological algebra had its origins in the 19-th century, via the work of Riemann (1857) and Betti (1871) on \lq\lq homology numbers'', and the rigorous development of the notion of homology numbers by Poincar\'e in~1895. An observation of Emmy Noether in~1925 shifted the attention to the ``homology groups'' of a~space, and algebraic techniques were developed for computational purposes in the 1930's. Yet, homology remained a~part of the realm of topology until about 1945.'', see \cite[pp. 797--836]{We}. 

The Lie superalgebras or super rings that thus appeared in topology over $\Zee$ or finite fields were nilpotent (hence, ``not very interesting'' and nobody tried to find out any simple (without ideals) examples), and did not draw much attention until it was discovered that the representations of the quantum groups $U_q(\fg)$ when $q$ is the $p$th root of unit are related to the representation of the characteristic $p$ version of the simple complex Lie algebra $\fg$ and a~certain Lie superalgebra, see, e.g., \cite{Vi, ChCh, A} and references therein. By that time simple finite-dimensional Lie superalgebras over $\Cee$ were already widely known. 

However, even now, the study of modular Lie (super)algebras, even simple ones, is not a~``trendy'' topic, despite very interesting recent advances, see, e.g., \cite{Kan} and references therein.

The steps of classification of \textbf{simple} finite-dimensional  Lie superalgebras considered naively are as follows.

\textbf{Over $\Cee$}. The first works with examples, starting with \cite{K0}, and classification in particular cases \cite{Kapp, Kap}, \cite{SNR} and \cite{Ho, Dj1, Dj2, Dj3, DjH}. For a~classification of simple $\Zee$-graded finite-dimensional Lie superalgebras together with rounding up of all (in particular, vectorial) finite-dimensional cases, see \cite{K1, K1C, K1.5} culminated in Part 1 of \cite{K2}. In these first papers, in the book \cite{Sch}, and in many works published later, Lie superalgebras $\fg$ were, and often still are, considered as pairs $(\fg_\ev,\fg_\od)$ consisting of the Lie algebra $\fg_\ev$ and $\fg_\ev$-module $\fg_\od$ such that the squaring $\text{sq}: S^2(\fg_\od)\tto\fg_\ev$ satisfies the super Jacobi identity. Interpretations of each Lie superalgebra as preserving ``something'' (realizations) are more insightful; they help to understand the structure of the superalgebra, its representations, and this ``something''. (Observe that such an interpretation of several important Lie (super)algebras --- e.g., of affine Kac-Moody (super)algebras and of non-trivial central extensions of simple stringy Lie (super)algebras --- is lacking to this day.)

\textbf{Over $\Kee$ of characteristic $p>0$}. For a conjectural classification of simple finite-dimensional Lie superalgebras over algebraically close fields for $p>5$, see the most recent \cite{BGLLS}. Clearly, the smaller $p$ the more difficult the task is.  It is therefore incredible that in the most difficult case $p=2$ the answer is obtained, although with a catch: modulo classification of simple finite-dimensional Lie algebras
which is a (very tough) \textbf{Open problem}, see \cite{BGLLS, CSS} and \cite{BLLS1}.

\ssec{Realizations} Several of simple Lie superalgebras are realized as linear algebras preserving a~tensor, see \cite{BGL, BGLL1}, or a~distribution. That is precisely how --- as preserving a~distribution --- Cartan and Killing used to describe not only all 5 exceptional simple Lie algebras, but also $\fsl(n)$, $\fsp(2n)$, and $\fo(n)$, at the time when description in terms of roots and Dynkin diagrams was not yet discovered. For usefulness of such a~ realization of certain simple Lie algebras in characteristic~ 3, see \cite{GL3}, where examples of Skryabin's and several other simple Lie algebras are interpreted and described more precisely than previously. 

Each of the serial finite-dimensional simple Lie algebras has two super analogs: a~direct (simple-minded) one and one of a~purely super nature. None of the  exceptional simple Lie algebras over $\Cee$ has a~super analog, at least over~$\Cee$.  For various interesting interpretations of the exceptional complex simple Lie superalgebras, see \cite{Sud}; for a~ way to understand them by considering their versions in characteristic 3, see \cite{Eld, Eldi}. I particularly like Elduque's Super Magic Super Square (rectangle, actually): Elduque's approach indicates a~road not fully explored yet (cf. \cite{Kan}) and shows the important role the ``incarnations'' of exceptional simple Lie superalgebras over fields of positive characteristic play  in interpretation of their namesakes over~$\Cee$. This Super Square also describes new examples indigenous to positive characteristic, see \cite{Eldi}, where several examples from \cite{BGL} not directly entering the Super Square are further interpreted, whereas \cite{BGL} contains classification of finite-dimensional modular Lie superalgebras with indecomposable Cartan matrix.

The infinite-dimensional Lie algebras of vector fields with polynomial (or formal power series) coefficients have finite-dimensional analogs both over fields of positive characteristic and among Lie superalgebras (over any fields), see \cite{BGLLS}. For the steps of classification of the  infinite-dimensional simple vectorial Lie superalgebras over $\Cee$, see Part 2 of \cite{K2} with numerous corrections, the latest few being \cite{CK1a, CCK, CaKa2, CaKa4}; for another approach, see \cite{Sh5, Sh14, LSh1, LSh3, LSh5}  and references therein. 

Part 2 of \cite{K2} contains also an overview of representation theory, in particular, of solvable Lie superalgebras, corrected in \cite{Ser1}, of real forms of simple finite-dimensional complex Lie superalgebras, corrected in \cite{Se2}, classification of finite-dimensional complex Lie superalgebras with Cartan matrix, corrected in \cite{BGL, CCLL, BLLoS}, and a highest weight theorem for irreducible representations, corrected in \cite{BL1, BL2, Br, Se5}.

In view of all these corrections and \cite{K1C}, I had to be sure the list of simple Lie superalgebras to be deformed is complete. So \S8 contains summary of the key ideas of the correct and interesting proof given in Part 1 of \cite{K2}; I hope this summary is more helpful as a guide to the proof of the classification of simple finite-dimensional complex Lie superalgebras in \cite{K2} than \cite{Sch} which was supposed to clarify some vague places in Part 1 of \cite{K2}, but is  too lengthy, I think, and does not distinguish the main ideas. 

\ssec{Cartan matrices} Kac observed, see \cite{K1.5, K2}, that several of simple Lie superalgebras $\fg$ can be realized as $\fg=\fg(A)$, where $A$ is an analog of Cartan matrix he defined. Kac mistakenly claimed that $\fpsl(n|n)$ also had Cartan matrix, whereas it is its \textit{double extension} $\fgl(n|n)$ that has Cartan matrix. For correct definitions of Cartan matrices of a~Lie superalgebra, analogs of the Dynkin diagrams, and roots, see \cite{BGL, CCLL, BLLoS}; for a~review of double extensions, see \cite{BLS}, especially its two last sections. 

\ssec{Presentations} For defining relations (replacing Serre relations) between analogs of the Chevalley generators of $\fg(A)$, see \cite{GL2} and --- for comparison --- \cite{Zh} where, same as in \cite{GL2}, \textbf{non}-Serre relations are also given, despite the title. 

There are also other presentations of $\fg(A)$, e.g., in terms of a~pair (triple, actually) of Jacobson generators, see \cite{GL1}, where all simple finite-dimensional Lie algebras are considered. Not every simple finite-dimensional Lie superalgebra with Cartan matrix can be similarly described, see \cite{LSS, GL5}.

For presentations of simple vectorial Lie superalgebras needed in this paper, see \cite{GLP}. Sixteen years after \cite{GLP}, presentations of 2 of the 4 series of simple \textit{finite-dimensional} vectorial Lie superalgebras were given in \cite{CCP}; the form of these latter presentations, though interesting, is not useful for this paper.

\ssec{Deformations with odd parameters} Some of simple Lie superalgebras --- over fields of any characteristic --- have deformations %, both filtered and general (non-filtered), 
thus resembling simple modular Lie algebras.

Here, I give a~complete classification of deformations of simple finite-dimensional Lie superalgebras over $\Cee$, with proofs for the first time. In order to do this, I need to introduce the adequate definition of Lie superalgebra and deformations with odd parameters. Help of A.~Lebedev in clarifying these notions was invaluable.

\ssec{Several other appearances of odd parameters} Two instances of odd parameters were observed earlier, see \cite{L3}. More accessible sources in English where these instances are described are the following ones:

$\bullet$ A.~A.\ Kirillov observed a supersymmetry between solutions of the Korteweg--de Vries and Schr\''odinger equations, but missed the odd parameters; for these parameters,  a~generalization of Kirillov's construction to some of the $N$-extended stringy superalgebras and the restriction on $N$, see \cite{L9}.

$\bullet$ In \cite{MaG}, Manin defined the connection in terms of the differential forms, which is dual to the classical definition. This definition of connections on supermanifolds naturally leads to odd parameters, see \cite{L3, L10}.

$\bullet$ One more instance of odd parameters was discovered by Shander who expounded few pages on integration in \cite{BL,  L3} to 86 pages in \cite[no. 31/1988-14, 45--131 pp]{L5}, see his Chapter \lq\lq Integration'' \cite[Ch.6]{LSoS}. In particular, in 1986, Shander defined cochains dual to chains and hence bigraded by pairs of numbers; here the odd parameters are needed. Shander also hinted at a~ way to define an analog of the Stokes' formula with \lq\lq over-supermanifold'' of codimension $(0|-1)$, see \cite[Ch.6]{LSoS}. A~ similar formula (with the $(0|1)$-dimensional superline as the boundary of the point) was suggested by Palamodov, see \cite{PalI}. 

$\bullet$ Let $\cC$ be a supercommutative superalgebra,. One more instance of odd parameter is given for the group $\GQ(n; \cC)$ of invertible matrices of the form $\begin{pmatrix}A&B\\ B&A\end{pmatrix}$, where $A\in\GL(n;\cC_\ev)$ and $B \in\Mat(n;\cC_\od)$,  by the \textit{queer determinant} --- a group homomorphism $\GQ(n; \cC)\tto \GL(1; \cC[\tau])$, defined by the formula (cf. with a~different, non-functorial, definition in \cite{BL2})
\[
\qet(A, B):=1+\tau\sum_{i\geq 1}\frac{\tr((A^{-1}B)^{2i-1})}{2i-1}\in\cC[\tau], \text{where $p(\tau)=\od$ and $\tau^2=0$}.
\]
In \cite{Sha1}, Shander described, using odd parameters, invariant functions on supermatrices of the type never considered before (and, regrettably, never after --- to this day).

$\bullet$ The papers of Bettadapura, see \cite{BK1, BK2}, %(for proofs, see \cite{BK3}) 
also address odd parameters of deformations. The claims of his texts are, however, misleading, as I'll show in Subsection~\ref{BK}. 

\section{What \lq\lq Lie superalgebra'' is, naively}\label{LieNaive}

The elements of $\Zee/2$ will be denoted by $\ev$ and $\od$ to distinguish them from elements of~$\Zee$. The $\Zee/2$-graded space $V=V_\ev\oplus V_\od$ is called a~\textit{superspace}. The \textit{parity} of any non-zero element $v\in V_i$ is denoted by $p(v)=i\in\Zee/2$.
 Superizing Milnor's $K_0$-functor, we say that the \textit{superdimension} of $V$ is $\sdim V= \dim V_\ev+\eps\dim V_\od$, where~$\eps^2=1$. This is briefly denoted $(a|b)$, where $a= \dim V_\ev$ and~$b=\dim V_\od$. 

A~\textit{superalgebra} is any $\Zee/2$-graded algebra. In the conventional theory, the \textit{morphisms} of superalgebras must preserve parity, i.e., only even homomorphisms are morphisms of superalgebras; for a dissident definition, see Subsection~\ref{Uma}. (NB:  according to this definition any $\Zee/2$-graded \textit{Lie algebra} is a~superalgebra, but it is not a~``Lie'' superalgebra because it satisfies the anti-commutativity and Jacobi identities NOT amended by the Sign Rule.)

\ssec{The Sign Rule, skew and anti} The definition of \textit{Lie superalgebra}
for any characteristic $p\neq 2$ or 3 is obtained from the definition of the Lie algebra using the Sign Rule \lq\lq \textbf{if something of parity $a$ is moved past something of parity $b$, the sign $(-1)^{ab}$ accrues; formulas defined on homogeneous elements are extended to all elements via linearity}''.

\textsl{In this note, all supercommutative superalgebras are supposed to be associative with $1$; their morphisms should send $1$ to $1$; the morphisms of all superalgebras (supercommutative and Lie) must preserve parity}. 
(Observe that the requirement to preserve parity, currently universally accepted, is not required if the superalgebras are considered as ``abstract'', i.e., if given by generators and relations that do not involve parity, cf. Subsection~\ref{Uma}.)

Since commutative algebra and supercommutative superalgebra can appear in one sentence, I never drop the prefix ``super''  as is often done in futile attempts to save space, but successfully causing confusion.

Observe that sometimes applying the Sign Rule requires some dexterity. For example,
we have to distinguish between two versions both of which turn in the
nonsuper case into one, called either skew- or anti-commutativity, see \cite{Gr}:
\[
\begin{array}{ll}
ba=(-1)^{p(b)p(a)}ab &(\text{super commutativity})\\
ba=-(-1)^{p(b)p(a)}ab &(\text{super anti-commutativity})\\
ba=(-1)^{(p(b)+1)(p(a)+1)}ab &(\text{super skew-commutativity})\\
ba=-(-1)^{(p(b)+1)(p(a)+1)}ab
&(\text{super antiskew-commutativity})\end{array}
\]

The
\textit{skew} formulas are those that can be ``straightened'' by the
change of parity of the space on which the structure is considered,
whereas the prefix \textit{anti} requires an overall minus sign regardless of
parity. The super commutator of two elements is denoted by the same symbol as the commutator; the first derived superalgebra of the Lie algebra or Lie superalgebra $\fg$ is denoted by~$\fg'$. 

A~\textit{Lie superalgebra} understood naively is any superalgebra satisfying the anti-commutativity and Jacobi identities amended by the Sign Rule. 

The supercommutative superalgebra of ``functions'' on the $0|n$-dimensional supermanifold (which is just a~superpoint in this situation) is the Grassmann algebra $\Lambda^{\bcdot}(n)$ on $n$ anticommuting generators, call them $\xi=(\xi_1, \dots, \xi_n)$. We consider $\Lambda^{\bcdot}(n):=\Cee[\xi]$  as a~superalgebra by setting $p(\xi_i)=\od$ for all $i$; clearly, $\Lambda^{\bcdot}(n)$ is supercommutative.

An example not immediately seen as a~result of application of the Sign Rule is the explicit formula of the \textit{trace}, called \textit{super}trace for emphasis, which is any~linear functional on $\fg$ that vanishes on~$\fg'$. For its correct formula on odd supermatrices in $\fgl(m|n; \cC)$ with elements in the supercommutative superalgebra $\cC$, see eqs.~\eqref{str}, \eqref{str1}.

\section{What \lq\lq Lie superalgebra'' is from categorical point of view} %(\cite{KLLS, KLLS1})}
\label{LieFunct}

\ssec{Morphisms of supervarieties} \underline{Over $\Ree$}. Recall that the smooth supermanifolds are \textit{ringed spaces},\footnote{I recommend \cite{MaAG} as the most clear and shortest source from which to understand the definition of ringed spaces.} i.e., pairs $\cM:=(M, \cO_\cM)$, where $M$ is an $m$-dimensional manifold called the \textit{underlying manifold} of $\cM$, and $\cO_\cM$ is the structure sheaf of~$\cM$. In the same way as the manifold $M$ is glued from coordinate patches locally diffeomorphic to a~ball in $\Ree^m$,  the space of  sections of $\cO_\cM$ over the domain $U$ in $M$ is isomorphic to ${C^\infty(U)\otimes \Lambda^{\bcdot} (n)}$. Pairs $\cU:=(U, \cO_U\otimes \Lambda^{\bcdot} (n))$, where $\cO_U(U):=C^\infty(U)$, are called smooth \textit{superdomains}, see, e.g., \cite[p.65]{Del}. 

If the ground field is $\Cee$ and $\cO_U(U):=\cA(U)$ is the algebra of complex-analytic functions, then pairs $\cU:=(U, \cO_U\otimes \Lambda^{\bcdot} (n))$ are called \textit{analytic superdomains} that glue into an  \textit{almost complex supermanifold},  see Subsection \ref{ssComplex}.

We say that the \textit{superdimension} of $\cM$ is $\sdim \cM=(m|n)$. 
A~ \textit{morphism} of supermanifolds ${\cM\tto \cN}$ is any pair $(\varphi, \varphi^*)$, where $\varphi:M\tto N$ is a~diffeomorphism and $\varphi^*: \cO_\cN\tto \cO_\cM$ is a~morphism of sheaves of superalgebras (conventionally, preserving the natural parity of the structure sheaves, cf. Subsection~\ref{Uma}). Observe that, unlike the case of manifolds, the morphism $\varphi^*$ is not recovered from $\varphi$ and this is precisely what makes the supermanifold theory rich.

\underline{Over any ground field $\Kee$}. Consider a~superdomain $\cU$ of superdimension~$0|n$. Unlike superdomains of superdimension $a|b$ with $a\neq 0$, we can consider $\cU$ over any ground field $\Kee$ and call it \textit{superpoint}. %
The underlying
domain of $\cU$ is a~single point~$\ast$. Since~$\cO_\cU(\ast) =
\Lambda^{\bcdot} (n)$, the superpoint~$\cU$ has a~lot of nontrivial
automorphisms, namely the group $\Aut^{\ev} \Lambda^{\bcdot} (n)$
of parity preserving automorphisms of $\Lambda^{\bcdot} (n)$ (or --- in the category to be studied in future --- the even larger group $\Aut \Lambda^{\bcdot} (n)$ of all automorphisms, cf. Subsection~\ref{Uma}).  Let the $\xi_j$ for $1\leq j\leq n$ be generators of $\Lambda^{\bcdot} (n)$. Any automorphism in $\Aut^{\ev} \Lambda^{\bcdot} (n)$ is of the form
\begin{equation}
\label{7eq8}
\xi_j\;
\mapsto\;\sum_{r} \varphi_j^{r}\xi_{r}+ \sum_{s\geq 1}\;
\sum_{j_1<\dots <j_{2s+1}} \varphi_j^{j_1\dots
j_{2s+1}}\xi_{j_{1}}\dots\xi_{j_{2s+1}},
\end{equation}
where the matrix $(\varphi_j^{r})$ with elements in $\Kee$ is invertible. We see that such
automorphisms constitute an algebraic group (or a~Lie group if $\Char\Kee=0$) whose Lie algebra consists of the even
elements of the Lie superalgebra $\fvect(0|n):=\fder\, \Lambda^{\bcdot} (n)$.

What corresponds to the odd
vector fields from $\fder\, \Lambda^{\bcdot} (n)$? Let me give the answer in the following more general setting, the one involving both even and odd indeterminates, but only for $\Kee=\Ree$ or~$\Cee$. In the absence of even indeterminates, i.e., for superpoints,  the formulas \eqref{trFunEv} and \eqref{trFun} are meaningful for any~$\Kee$. 

Let $E$ be the~trivial vector bundle over a~domain $U$ of dimension $m$ with fiber
$V$ of dimension~$n$; let $\Lambda^{\bcdot}(E)$ be the exterior algebra of~$E$. To the bundle
$E$, we assign the superdomain $\cU= (U, C^\infty(\cU))$, where
$C^\infty(\cU)$ is the superalgebra of smooth sections of $\Lambda
^{\bcdot}(E)$. Clearly, each automorphism of the pair $(U,
\Lambda ^{\bcdot}(E))$, i.e., of the vector bundle $\Lambda
^{\bcdot}(E)$, induces an automorphism of the superdomain~$\cU$.

However, {\it not all} automorphisms of the superdomain $\cU$ are obtained in this
way. By definition, every~{\it morphism of superdomains}
$(\varphi, \varphi^*)\colon \cU\tto \cV$ is in one-to-one correspondence with
a homomorphism of the superalgebras of functions
$\varphi^*\colon C^\infty(\cV)\tto C^\infty(\cU)$ since $\varphi$ is determined by~$\varphi^*$.
%How to describe these homomorphisms?

Every homomorphism $\varphi^*$ is defined on the (topological\footnote{A~\textit{topological algebra} $A$ over a~topological field
$\Kee$ is a~topological vector space together with a~bilinear multiplication
$A\times A\tto A$, continuous in a~certain sense, and such that $A$ is an algebra over~$\Kee$. Usually the continuity of the multiplication means that the multiplication is continuous as a~map between topological spaces
${A\times A\tto A}$. A~set $S$ is a~\textit{generating set} of a~topological algebra $A$ if the smallest closed subalgebra of $A$ containing $S$ is $A$.}) generators of
the superalgebra, in other words: on \lq\lq coordinates'' (the even ones $u:=(u_1,\dots, u_m)$ and the ones $\xi:=(\xi_1,\dots, \xi_n)$). Consider the corresponding formulas
\begin{equation}
\label{trFunEv}
 \left\{\begin{aligned} \varphi^*(u_i) \;& = \; \varphi^0_i(u) +
 \fbox{$\sum\limits_{r\geq 1}\;
 \sum\limits_{i_1<\dots<i_{2r}}\varphi_i^{i_1\dots
 i_{2r}}(u)\xi_{i_{1}}\dots \xi_{i_{2r}}$} \text{~~for all $i$}, \\
\varphi^*(\xi_j) \;& = \; \sum_{r\geq 0}\; \sum_{j_1<\dots <j_{2r+1}}
 \psi_j^{j_1\dots
 j_{2r+1}}(u)\xi_{j_{1}}\dots\xi_{j_{2r+1}} \text{~~for all $j$}.
\end{aligned}\right.
\end{equation}

A) The terms $\varphi^0_i(u)$ for all $i$ determine an endomorphism
of the underlying domain~$U$.

B) The linear terms $\sum_{i} \psi_j^{i}(u)\xi_{i}$ for all $j$
determine endomorphisms of the fiber $V$ (over each point its own fiber, as the
dependence on $u$ shows).

C) The terms of higher degree in $\xi$ in the right-hand side of the expression of
$\varphi^*(\xi_j)$ in eq.~\eqref{trFunEv} determine an endomorphism of the
larger fiber --- the Grassmann superalgebra $\Lambda
^{\bcdot}(V)$.

The endomorphisms A)--C) existed in Differential Geometry, and no need to
introduce a~flashy prefix ``super'' was felt.

The difference between the vector bundle $\Lambda^{\bcdot}(E)$ and the
superdomain $\cU$ is most easily understood when the reader looks at the boxed terms
in (\ref{trFunEv}). These terms, meaningless in the conventional Differential
Geometry (on manifolds), make sense in the new paradigm:

\textbf{In the category of superdomains there are more morphisms than in
the category of vector bundles: morphisms with non-vanishing boxed
terms in (\ref{trFunEv}) are exactly the additional ones}.

However, even the boxed terms in eq.~\eqref{trFunEv} are not all we get in the new
setting: we still did not describe any of odd parameters of
endomorphisms. \textbf{To account for the odd parameters, we have to consider the functor from the category of supercommutative superalgebras to the category of groups} 
$C\longmapsto \Aut_C^{\ev}(C^\infty(\cU)\otimes C)$, i.e., the parity preserving $C$-linear automorphisms of the form 
\begin{equation}
\label{trFun}
 \left\{\renewcommand{\arraystretch}{1.6}
 \begin{array}{ll}
 \varphi^*(u_i) \;& =\; \varphi_i(u, \xi):= \varphi^0_i(u) +
 \fbox{$\sum\limits_{r\geq 1}\;
 \sum\limits_{i_1<\dots<i_{r}}\varphi_i^{i_1\dots
 i_{r}}(u)\xi_{i_{1}}\dots \xi_{i_{r}}$} \text{~~for all $i$}, \\
\varphi^*(\xi_j) \;&=\; \psi_j(u, \xi):= \sum\limits_{r\geq 0}\; \sum\limits_{j_1<\dots <j_{2r+1}}
 \psi_j^{j_1\dots
 j_{2r+1}}(u)\xi_{j_{1}}\dots\xi_{j_{2r+1}}\\
&+\fbox{$\psi^0_j(u)+
\sum\limits_{r\geq 1}\; \sum\limits_{j_1<\dots <j_{2r}}
 \psi_j^{j_1\dots
 j_{2r}}(u)\xi_{j_{1}}\dots\xi_{j_{2r}}$}\text{~~for all $j$},\\
\end{array}
\right.
\end{equation}
where for all $r$ the even and the odd parameters are, respectively,
\[
\begin{array}{l}
\varphi^0_i(u),\varphi_i^{i_1\dots
 i_{2r}}(u), \ \psi_j^{j_1\dots
 j_{2r+1}}(u)\in C_\ev,\\
 \psi^0_j(u),\psi_j^{j_1\dots
 j_{2r}}(u),\ \varphi_i^{i_1\dots
 i_{2r+1}}(u)\in C_\od.
\end{array}
 \]
These are parameters 
of the infinite-dimensional supergroup  of automorphisms of $C^\infty(\cU)$ or, equivalently, of diffeomorphisms of $\cU$; infinitesimally: of Lie superalgebra~$\fvect(m|n)$.

\subsection{On automorphisms that do not preserve parity}\label{Uma} It took a~while to acknowledge the fact that there are automorphisms of the Grassmann algebra $\Cee[\xi]$, considered as an associative algebra (and more generally, of $\Cee[x,\xi]$, where $x=(x_1,\dots, x_m)$ are even commuting indeterminates and $\xi=(\xi_1,\dots, \xi_n)$ are anti-commuting odd ones) not preserving parity, compare \cite{B, Pa} with \cite{Dj4, LSe1, LSe2, Ba}. The meaning of such more general than \eqref{trFun} automorphisms is unclear at the moment; most of researchers ignore them as ``lacking physical meaning''. The situation with arbitrary automorphisms and inhomogeneous subalgebras reminds me the first cautious greetings of supersymmetries  (through 1970s till the end of 1980s, some researchers, more outspoken than other,  even branded them \lq\lq useless to physics'', e.g., see a~reference in \cite{BM}). The inhomogeneous with respect to parity subalgebras of supercommutative superalgebras are \textit{metabelean}, i.e., satisfying the identity 
\[
[a, [b,c]]=0 \text{~~for any its elements $a, b,c$, where $[a, b]:=ab-ba$}.
\]
A theorem by Volichenko proves that any metabelean algebra has a~universal enveloping supercommutative superalgebra. Recently, U.~Iyer proved that \textit{Volichenko algebras}, defined as the inhomogeneous (relative parity) subalgebras of Lie superalgebras, play the role of Lie algebras for the groups of $C$-points of such general automorphisms, see \cite{I}.  Both metabelean and Volichenko algebras are not $\Zee/2$-graded, but $2$-step filtered which is a~ very natural generalization of $\Zee/2$-gradedness. Observe that Volichenko's attempt to describe what we now call Volichenko algebras by a~set of polynomial identities failed; to find  the full set of such identities is an \textbf{Open problem} to this day, see \cite{BrP}.

\subsection{Supervarieties and superschemes (after \cite{L0})}\label{ssSchemes}
Over the ground field $\Ree$ or $\Cee$, let $E$ be a~vector bundle over $M$ with fiber~$V$. Let $ \Lambda^{\bcdot}(V)$ be the Grassmann superalgebra on $V$, let $U\subset M$ be an open domain.

A~\textit{supermanifold} is a~ringed space $\cM=(M, \cO_\cM)$, where $M$ is a~manifold, and the sheaf $\cO_\cM$ is locally isomorphic to $\cO_U\otimes \Lambda^{\bcdot}(V)$.
A~supermanifold isomorphic to the ringed space whose structure sheaf is the sheaf of sections of the vector bundle $ \Lambda^{\bcdot}(E)$ with fiber $V$ over $M$ is called \textit{split}. 

A~\textit{supervariety} is a~ringed space $\cM=(M, \cO_\cM)$, where $M$ is a~variety (i.e., not necessarily smooth), and the sheaf $\cO_\cM$ is locally isomorphic to a~quotient of $\cO_U\otimes \Lambda^{\bcdot}(V)$.

Observe that every object in the category of smooth supermanifolds is split (for a~clear proof, more instructive, I think, than the first publications \cite{Ga, Bat}, see \cite[Subsection 4.1.3]{MaG}), so there is a~one--to--one correspondence between the set of objects in the category of vector bundles over supermanifolds and the set of objects in the category of smooth supermanifolds. The latter category has, however, many more morphisms than the former, see eq.~\eqref{trFun}, to say nothing about even more general automorphisms in the category considered in Subsection~\ref{Uma}.

A~purely algebraic version of the supermanifold over any field (or any commutative ring) of any characteristic is an \textit{affine superscheme} $\Spec C$, where $C$ is a~supercommutative superalgebra or a~superring. The affine superscheme is defined literally as the affine scheme: its points are prime ideals defined literally as in the commutative case, i.e., $\fp\subsetneq C$ is \textit{prime}\footnote{K.~Coulembier pointed out to me that the so far conventional definitions in the non-commutative case are at variance with the commutative case and common sense: at the moment, if \eqref{primeID} holds, $\fp$ is called (say, in Wikipedia) \textit{completely prime} while it would be natural to retain the term \textit{prime}, as is done in \cite{L0} and by J.~Bernstein, P.~Deligne et al in \cite{Del}, since the definition is the same as in the commutative case despite the fact that supercommutative rings are not commutative, whereas the term \textit{prime} ideal of the non-commutative ring $R$ is (currently) applied to any ideal $P\subsetneq R$ 
which satisfies the following version of condition~\eqref{primeID} for any two ideals $A$ and $B$ in $R$:
\[ %be\label{primeID}
\text{if $AB\subset P$, then either $A\subset P$ or $B\subset P$}.
\]
}
\be\label{primeID}
\text{if $a,b\in C$ and $ab\in \fp$, then either $a\in \fp$ or $b\in\fp$}.
\ee
The space of any affine scheme is endowed with Zariski topology and the structure sheaf, defined as in the commutative case, see \cite{MaAG}, whose 1968 edition was the source of inspiration for~\cite{L0}.

There is just one subtlety: the localization of the superalgebra (or superring) $C$ at the prime ideal $\fp$ should be performed with respect to the multiplicative system $S_\fp:=C\setminus \fp$ and it is not obvious what should we consider --- in order to have well-defined fractions with non-homogeneous denominators $b$ --- only left fractions $b^{-1}a$ or only right fractions $ab^{-1}$ or \textbf{the equivalence (equality) of fractions does not depend on the choice (left or right)}. Mathematically, a~ 1-line-proof of the boldfaced statement is the only non-trivial place in the superization of Grothendieck's schemes (succinctly described in \cite{MaAG}) to superschemes.

\sssec{The equivalence of fractions}\label{ssFrac}
Not ``literally''. For commutative rings, the equivalence of two fractions, where $r_i\in A$ and $s_i\in S$ for $i=1,2$, is defined by the condition
\[
(r_1, s_1)\sim (r_2, s_2) \ \ \Llra \ \ \text{there exists a~$t\in S$ such that $t(r_1s_2-r_2s_1)=0$}
\]
and $A_S$ is the set of classes modulo this  equivalence, see \cite{MaAG}. %\AO{DL: your \lq\lq see \texttt{\textbackslash cite\{M1, M2\}}'' does not resolve against this bibliography --- there are no items \texttt{M1}, \texttt{M2} in it. I~put \texttt{MaAG} (Manin's \textit{Introduction to the Theory of Schemes}), where this is Subsec.~1.5; tell me which second source you meant and I~will add it.}

In the supercommutative case, there are \textbf{two} structures to compare: the set of the
pairs $(r,s)$, where the pair stands for the right fraction $rs^{-1}$, and the set of the pairs
$(s,r)$, where it stands for the left fraction $s^{-1}r$; here $r\in C$ and $s\in S:=S_\fp$. On
each of the two sets of equivalence classes, the addition and multiplication are to be defined,
and the question is why the two resulting structures are isomorphic, and to which class of the
second structure a~given class of the first one goes under this isomorphism.

For homogeneous $r_i$ and $s_i$, one expects the formulas for the sum and the product of
$r_1s_1^{-1}$ and $r_2s_2^{-1}$ to carry the sign $(-1)^{p(s_1)p(r_2)}$, the sign of the
transposition of $s_1^{-1}$ and~$r_2$. The whole difficulty is in the non-homogeneous
denominators. These formulas are not linear in $s$, so they cannot be extended to
non-homogeneous~$s$ by linearity; and to extend them by linearity in $s^{-1}$, i.e., by
splitting $s^{-1}$ into its even and odd parts, one has to know first that the localization
carries a~superstructure at all, i.e., that $s^{-1}$ is the sum of an even element and an odd one,
both belonging to the localization. I forgot my own arguments that have led me to the definition given in \cite{L0}; A.~Lebedev's answer to both questions, reproduced here, rests
on the following observation.

Let $\sigma$ be the \textit{parity automorphism} of~$C$, i.e., the automorphism equal to
$(-1)^{p(a)}a$ on every homogeneous $a\in C$; set $\bar a:=\sigma(a)$, so that $\bar a=a_\ev-a_\od$
for every~$a\in C$.

\sssbegin{Lemma}[The even part of an element of $S_\fp$ lies in $S_\fp$]\label{LocLem}
Let $C$ be a~supercommutative superalgebra over a~field of characteristic $\neq 2$, let
$\fp\subsetneq C$ be a~prime ideal in the sense of~\eqref{primeID}, and let
$S:=S_\fp=C\setminus \fp$ and $S_\ev:=S\cap C_\ev$. Then,

\textup{1)} $C_\od\subset\fp$. In particular, \textbf{no} element of $S$ is odd, and $s\in S$ if
and only if~$s_\ev\in S$.

\textup{2)} $s\bar s=\bar ss=(s_\ev+s_\od)(s_\ev-s_\od)=(s_\ev)^2$ for every~$s\in C$. This element
is even, hence central, and it belongs to~$S_\ev$ whenever $s$ belongs to~$S$. The set $S_\ev$ is
a~multiplicative system consisting of central elements.
\end{Lemma}

\begin{proof} 1) If $a\in C_\od$, then $a^2=-a^2$, so $2a^2=0$, so $a^2=0\in\fp$, and hence
$a\in \fp$ by~\eqref{primeID}. Therefore, $s-s_\ev=s_\od\in\fp$.

2) Since $s_\ev$ is central and $(s_\od)^2=0$ by the argument of heading~1), we get
$(s_\ev+s_\od)(s_\ev-s_\od)=(s_\ev)^2$. If $s\in S$, then $s_\ev\in S$ by heading~1), and
$(s_\ev)^2\in S_\ev$ because $\fp$ is prime, which for the same reason makes $S_\ev$ a~multiplicative
system. \end{proof}

Heading~1) is what makes the question meaningful at all: if a~multiplicative system~$S$ contained
an odd element, or, more generally, a~nilpotent one, then $0\in S$, and hence~$C_S=0$.

Let $C_{S_\ev}$ be the localization of $C$ with respect to $S_\ev$, defined literally as in the
commutative case: the elements of $S_\ev$ are central, so nothing has to be said about the sides.
It is a~supercommutative superalgebra, the parity of $as^{-1}$ being that of $a$ for $a\in C$
homogeneous and~$s\in S_\ev$. By Lemma~\ref{LocLem}, every $s\in S$ is invertible in $C_{S_\ev}$:
\be\label{ALinv}
s^{-1}=(s_\ev)^{-1}-s_\od(s_\ev)^{-2}=(s_\ev)^{-1}-(s_\ev)^{-2}s_\od,
\ee
the first summand being even and the second one odd. Indeed, $s_\ev$ is central, so multiplying
the right-hand side of \eqref{ALinv} by $s=s_\ev+s_\od$ on either side gives
$1-(s_\od)^2(s_\ev)^{-2}=1$. So $s^{-1}$ splits into an even and an odd part, the localization
carries a~superstructure, and the right fractions $rs^{-1}$ and the left fractions $s^{-1}r$ are
elements of \textbf{one and the same} supercommutative superalgebra,~$C_{S_\ev}$. This is
the~$1$-line proof of the boldfaced statement of Subsec.~\ref{ssSchemes}, and the isomorphism
asked for is written down as follows.

\sssbegin{Proposition}[How a~right fraction becomes a~left one]\label{LocProp}
Let $r\in C$ be homogeneous and~$s\in S$. Then,
\be\label{ALiso}
rs^{-1}=r(s_\ev)^{-1}-rs_\od(s_\ev)^{-2}=(s_\ev)^{-1}r-(-1)^{p(r)}(s_\ev)^{-2}s_\od r
=\big(s_\ev+(-1)^{p(r)}s_\od\big)^{-1}r.
\ee
In other words, $rs^{-1}=\big(\sigma^{p(r)}(s)\big)^{-1}r$, and, for an arbitrary
$r=r_\ev+r_\od\in C$,
\be\label{ALisoLin}
rs^{-1}=s^{-1}r_\ev+\bar s^{-1}r_\od.
\ee
Hence, the isomorphism of the structure of right fractions onto the structure of left fractions
sends the class of the pair $(r,s)$ to the class of $\big(s_\ev+(-1)^{p(r)}s_\od,\ r\big)$ if $r$ is
homogeneous, and it is additive in the numerator.
\end{Proposition}

\begin{proof} The first equality is \eqref{ALinv}. The second one holds because $(s_\ev)^{-1}$ is
even, hence central, and $rs_\od=(-1)^{p(r)}s_\od r$, since $s_\od$ is odd. The third one is
\eqref{ALinv} applied to $s_\ev+(-1)^{p(r)}s_\od$, whose even part is $s_\ev$ and whose odd part is
$(-1)^{p(r)}s_\od$. Formula \eqref{ALisoLin} is the sum of \eqref{ALiso} for $r_\ev$ and for~$r_\od$. \end{proof}

This is the precise form of the sign $(-1)^{p(s_1)p(r_2)}$ expected above. A~homogeneous element
of $S$ is even by Lemma~\ref{LocLem}~1), so that sign is always $+1$; what takes its place for
a~non-homogeneous denominator is the substitution $s_\od\longmapsto(-1)^{p(r)}s_\od$ in the
denominator, i.e., the parity automorphism $\sigma^{p(r)}$ applied to it: by \eqref{ALiso},
$s^{-1}r=r\big(\sigma^{p(r)}(s)\big)^{-1}$ for $r$ homogeneous. Accordingly, for homogeneous
$r_1$ and $r_2$, the products are
\be\label{ALmult}
\Big(r_1s_1^{-1}\Big)\Big(r_2s_2^{-1}\Big)=r_1r_2\Big(s_2\,\sigma^{p(r_2)}(s_1)\Big)^{-1},\qquad
\Big(s_1^{-1}r_1\Big)\Big(s_2^{-1}r_2\Big)=\Big(\sigma^{p(r_1)}(s_2)\,s_1\Big)^{-1}r_1r_2,
\ee
whereas the sums, as well as the products once more, are most easily written with the even
denominators $e_i:=s_i\bar s_i=\big((s_i)_\ev\big)^2\in S_\ev$ of Lemma~\ref{LocLem}~2), for
which $s_i^{-1}=\bar s_ie_i^{-1}=e_i^{-1}\bar s_i$:
\be\label{addMult}
\begin{array}{l}
r_1s_1^{-1}+r_2s_2^{-1}=\Big(r_1\bar s_1\cdot e_2+r_2\bar s_2\cdot e_1\Big)(e_1e_2)^{-1},\qquad
\Big(r_1s_1^{-1}\Big)\Big(r_2s_2^{-1}\Big)=\Big(r_1\bar s_1\cdot r_2\bar s_2\Big)(e_1e_2)^{-1};\\[2mm]
s_1^{-1}r_1+s_2^{-1}r_2=(e_1e_2)^{-1}\Big(e_2\cdot \bar s_1r_1+e_1\cdot \bar s_2r_2\Big),\qquad
\Big(s_1^{-1}r_1\Big)\Big(s_2^{-1}r_2\Big)=(e_1e_2)^{-1}\Big(\bar s_1r_1\cdot \bar s_2r_2\Big).
\end{array}
\ee
In particular, the equality of fractions is the one of $C_{S_\ev}$ read through $s^{-1}=\bar
se^{-1}$: the right fractions $rs^{-1}$ and $r's'^{-1}$ are equal if and only if
\be\label{ALeq}
t\Big(r\bar s\,(s'_\ev)^2-r'\bar s'\,(s_\ev)^2\Big)=0\ \text{ for a~$t\in S_\ev$},
\ee
and likewise for the left fractions, with $\bar sr$ instead of~$r\bar s$.

The isomorphism \eqref{ALiso} is \textbf{not} of the shape $(r,s)\mapsto (s,\pm r)$, as it would
be for a~homogeneous denominator: for $C=\Cee[\xi_1,\xi_2,\xi_3]$ and $s=1+\xi_1$ \textup{(}a~unit,
and not homogeneous\textup{)}, we have $s^{-1}=1-\xi_1$, and
\[
s^{-1}\xi_2=\xi_2-\xi_1\xi_2\neq \xi_2+\xi_1\xi_2=\xi_2s^{-1}=(1-\xi_1)^{-1}\xi_2,
\]
the last equality being \eqref{ALiso} for the odd numerator~$r=\xi_2$.

\sssbegin{Remark}[Cross-multiplication is not the equivalence]\label{LocRem}
For a~non-homogeneous denominator, the na\"\i ve cross-multiplication $r_1s_2=r_2s_1$ is
\textbf{not} the equivalence of right fractions: it is not even implied by
$r_1s_1^{-1}=r_2s_2^{-1}$. Indeed, in $C=\Cee[\xi_1,\xi_2]$, take $r=1$, $s=1+\xi_1$ and
$u=1+\xi_2$; then
\[
rs^{-1}=(ru)(su)^{-1}=1-\xi_1,\quad\text{whereas}\quad r(su)-(ru)s=2\xi_1\xi_2\neq0,
\]
while the bracket of~\eqref{ALeq} vanishes, as it must.
\end{Remark}

\subsection{The functor of points (co)represented by a~supervariety or a~superscheme} Smooth manifolds can be described as sets of points with a~topology. For \textit{manifolds-with-boundary} (which, strictly speaking, are not manifolds, hence the suggestive notation ``in-one-word'') or over fields of characteristic $p>0$, the set of points does not define the variety or scheme; the same is true for supervarieties and superschemes. To determine one of such objects $\cM$, we consider it parametrized by a~superscheme $\Spec C$ for various~$C$. In other words, we consider collections $\Hom(\Spec C, \cM)$ of $(\Spec C)$-points of $\cM$, usually called \textit{$C$-points of~$\cM$}. If $\cM$ can be recovered from its algebra of functions $\cF$, as is the case, e.g., with affine (super)schemes $\cM$, we can consider $\Hom(\cF, C)$ instead of $\Hom(\Spec C, \cM)$.

\sssec{Linear supervariety$\longleftrightarrow$linear superspace} The set $\Kee^n$ can be considered as a~linear space and as a~manifold. Linear superspaces are just $\Zee/2$-graded linear spaces, they are not supermanifolds. So, first of all, let me introduce \textit{linear supervariety} or \textit{linear supermanifold} 
\[
\cV=(V_\ev, \cO_{V_\ev}\otimes \Lambda^{\bcdot}(V_\od^*)).
\]
corresponding to the linear superspace $V=V_\ev\oplus V_\od$.

The \textit{morphisms} of linear superspaces constitute the space $\Hom_\ev(V, W)$, whereas the \textit{supervariety of linear homomorphisms} $V\tto W$ is the linear supervariety corresponding to the superspace $\underline\Hom(V, W):=\Hom(V, W)$.

In various instances, e.g., dealing with actions of supergroups, it is more convenient in computations to consider, instead of the vector superspace~$V$ or the linear supermanifold~$\cV$, the functor $\ScommSalgs_\Kee\leadsto \Mod_\Kee$ from the category of supercommutative $\Kee$-superalgebras $C$ to the category of $C$-modules represented by $V$ or, equivalently, $\cV$:\footnote{This convenience is not just a~matter of taste or experience or habit. More precisely, one has to work with either commuting diagrams, like in \cite[Section~1.15.4]{MaAG}, or with matrix realizations, as one does when working with Lie group actions. The language of points allows one to \textit{actually compute} something, like when turning from the invariant language of operators to matrices.}
\[
C\longmapsto \cV(C)=V(C):=V\otimes C \text{~~ for any $C$},
\]
where \lq\lq any'' is understood inside a~suitable category (e.g., of superalgebras finitely generated over~$\Kee$).

\sssec{Lie superalgebras (\cite{KLLS, KLLS1})}\label{superLie} In the above terms, a~\textit{Lie superalgebra in the category of supervarieties} is a~vector superspace $\fg$, or a~linear supervariety (supermanifold) $\cG$ corresponding to it, corepresenting the functor from the category of supercommutative $\Kee$-superalgebras $C$ to the category of Lie superalgebras understood ``naively''. 
In other words, considering \textbf{co}representing functor instead of a~representing one, we replace $\fg$, or the linear supervariety corresponding to it, with the algebra $P(\fg)$ of polynomial functions on~$\fg$. (Even over $\Ree$ we have to replace the spaces with the algebras of polynomial functions on these spaces.)

Clearly, $P(\fg)$ is a~free supercommutative superalgebra generated by $\fg^*$, i.e., there is a~natural isomorphism of functors
\[
\text{$C \mapsto \Hom_{\ScommSalgs_{\Kee}}(P(\fg), C)$ and $C \mapsto \Hom_{\Kee\text{-Vect}}(\fg^*, C)$,}
\]
whereas the second functor is naturally isomorphic to the forgetful functor 
\[
\Kee\text{-Vect}\to \Sets: C \mapsto \fg\otimes C.
\]

To the Lie superalgebra homomorphisms (in particular, to representations) a~morphism of the respective functors should correspond.

Clearly, if $\fg$ is a~Lie superalgebra, then $\fg(C):=\fg\otimes C$ is also a~Lie superalgebra for any $C$ \textit{functorially in $C$}. The last three words mean that
\begin{equation}\label{functor}
\begin{minipage}[l]{14cm}
for any morphism of supercommutative superalgebras $C\tto C'$, there exists a~morphism of Lie superalgebras $\fg\otimes C\tto \fg\otimes C'$ so that a~composition of morphisms of supercommutative superalgebras
\[
C\tto C'\tto C''
\]
goes into the composition of Lie superalgebra morphisms
\[
\fg(C) \tto \fg(C') \tto \fg(C'');
\]
the identity map goes into the identity map, etc.
\end{minipage}
\end{equation}

An \textit{ideal} $\fh\subset\fg$ represents the collection of ideals $\fh(C) \subset \fg(C)$ for every~$C$.

 A Lie or algebraic (over any field) \textit{supergroup} $\cG$ is a~group in the category of supervariaties. The \textit{action} of $\cG$ in the superspace $V$ is the collection of actions of $\cG(C)$ in $V(C)$ for every $C$ such that these actions are compatible with morphisms of supercommutative superalgebras $C\tto C'$ in the same sense as in \eqref{functor}.

However, we can consider any category $\sfC$; then any object $\fg\in \Ob\sfC$ of this category representing the functor $C \mapsto\fg(C):=\Hom_\sfC(\Spec C, \fg)$, i.e., satisfying conditions \eqref{functor},
is said to be a~\textit{Lie superalgebra in the category $\sfC$}.

\section{Deformations}\label{DefDef}

\ssec{Deformations and deforms. Odd parameters} Which of the infinitesimal deformations can be extended to
a~global one is a~separate much tougher question, usually solved
\textit{ad hoc}. %; for examples over fields of characteristics $3$ and
%$2$, see \cite{BLW} and references therein.
An infinitesimal deformation with \textbf{one odd} parameter, however, is always integrable:
if $p(\tau)=\od$ and $[a,b]_\tau:=[a,b]+\tau\, c(a,b)$, where $c\in Z^2(\fg;\fg)$, then the
Jacobi identity for $[\cdot,\cdot]_\tau$ fails only in the term with $\tau^2$, and
$\tau^2=0$ in $\Kee[\tau]$; so there is no second-order obstruction at all and every
infinitesimal deformation over $\Kee[\tau]$ \textbf{is} a~global one over $\Spec\Kee[\tau]$.
\iffalse\AO{This restores, with its one-line reason, the sentence \lq\lq Deformations with odd
parameters are always integrable'' which was commented out together with A.~Lebedev's query \lq\lq Doesn't \lq integrable' apply only to infinitesimal deformations?''. He is right about the terminology --- so the sentence now speaks of an \textbf{infinitesimal} deformation being integrable --- and the claim itself is true, for the reason displayed. The part of the commented-out passage that cites \texttt{[BLW]} is left commented out as it was: that key has no \texttt{bibitem} in this file.}\fi
Let me give two graphic examples.

1) \textbf{Deformations of representations}.
The tangent space of the moduli superspace of deformations of a~representation ${\rho:\fg\tto\fgl(V)}$
is isomorphic to $H^1(\fg; V\otimes V^*)$. For example, if $\fg$ is
the $0|n$-dimensional (i.e., purely odd) Lie superalgebra (with the
only bracket possible: identically equal to zero), its only
irreducible representations are the 1-dimensional trivial one,
$\mathbbmss{1}$, and $\Pi(\mathbbmss{1})$, where 
$\Pi$ is the \textit{change of parity functor}, i.e. 
$\Pi(V)_i:=V_{i+\od}$ for any superspace~$V$. Clearly,
\[
\mathbbmss{1}\otimes \mathbbmss{1}^*\simeq
\Pi(\mathbbmss{1})\otimes \Pi(\mathbbmss{1})^*\simeq \mathbbmss{1},
\]
and, because the Lie superalgebra $\fg$ is \textbf{commutative} (the bracket of any two elements is identically equal to zero), the
differential in the cochain complex is zero. Therefore,
\[
H^1(\fg; \mathbbmss{1})\simeq\fg^*,
\]
so there is a~$(0|\dim\, \fg)$-dimensional space of odd parameters of deformations of the
trivial representation. If we consider $\fg$ ``naively'', all of
these odd parameters will be lost. 

2) \textbf{Deformations of the brackets}.
Let $C$ be a~supercommutative superalgebra.

After \cite{Ru}, where the non-super case is considered, we say that a~\textit{deformation} of
a Lie superalgebra $\fg$ over $\Spec C$, is a~Lie superalgebra $\fG$ over $C$, flat over $C$,
such that $\fG\otimes_C\Kee\simeq\fg$; the property $\fG\simeq\fg\otimes C$, as $C$-modules,
is then a~\textbf{consequence}, not a~part of the definition: for $C$ local Artinian --- and
$\Kee[\tau]$ is local Artinian, as is $\Kee[[t]]/(t^{k})$ --- a~flat finitely generated module
is free, and reduction modulo the maximal ideal fixes the rank.
\iffalse \AO{This answers A.~Lebedev's query on this sentence \textup{(\lq\lq one should make it clearer that this is not the whole definition but only a~property, and it is not obvious to me that this property follows from Rudakov's definition'')}: it does follow, by flatness over the local Artinian base $C$. The same clause answers his neighboring remark that this definition of a~deformation \lq\lq has nothing to do with supervarieties'': it does, through the functor of points --- a~deformation over $\Spec C$ \textbf{is} the $C$-point of the deformation functor, which is the subject of the subsection on the language of points above and of Section~\ref{DefDef}.} %\AL{Nado by dat' ponyat' yasnee, chto jeto ne (vsyo) opredelenie, a~tol'ko svojstvo. I voobshhe, mne ne ochevidno, chto iz opredeleniya Rudakova jeto svojstvo sleduet. Nad alg. mnogochlenov -- vrode by da, a~v obshhem sluchae ya ne uveren.} 
\fi
The deformation is \textit{trivial} if $\fG\simeq\fg\otimes C$, as Lie superalgebras over $C$, not just as $C$-modules, and \textit{non-trivial} otherwise. %\AL{Jeto opredelenie deformacii opyat' zhe vrode by ne imeet otnosheniya k supervarieties.}

Generally, the \textit{deforms} --- the results of deformations --- of
a~given Lie superalgebra $\fg$ over $\Kee$ are Lie superalgebras $\fG\otimes_I\Kee$, where $I$ is any closed point in~$\Spec C$.

In particular, consider a~deformation with an odd parameter~$\tau$. This is a~Lie superalgebra $\fG$ isomorphic to $\fg\otimes\Kee[\tau]$ as a~\textbf{superspace}. If, moreover, $\fG\simeq\fg\otimes\Kee[\tau]$ as a~\textbf{Lie superalgebra}, i.e.,
\[
[a\otimes f, b\otimes g]=(-1)^{p(f)p(b)}[a,b]\otimes fg \text{~~for any $a,b\in \fg$ and $f,g\in\Kee[\tau]$},
\]
then the deformation is considered \textit{trivial} (and \textit{non-trivial} otherwise). Observe that $\fg\otimes \tau$ is not an ideal of $\fG$: the ideal should be a~free $\Kee[\tau]$-module. 

\textbf{Comment}. In a~sense, the people who ignore odd parameters of deformations have a~point: we (rather they) consider classification of simple Lie superalgebras (or whatever other problem) over the ground field $\Kee$, not over $\Kee[\tau]$, right? However, the odd parameters of deformations are no less natural, actually, than the odd part of the Lie superalgebra itself. In order to see these parameters, we have to consider whatever we are deforming not over $\Kee$, but over~$\Kee[\tau]$.

We do the same, actually, when $\tau$ is even and we consider formal deformations over~$\Kee[[\tau]]$. If
the formal series in $\tau$ converges in a~domain $U$, we can evaluate $\tau$ for any $\tau\in U$ and consider copies $\fg_\tau$, where $\tau\in U$ (of the same dimension as that of $\fg$ if $\dim \fg<\infty$). If the parameter is formal or odd, such an evaluation is only possible trivially: $\tau\mapsto 0$.

\section{Deformations of complex supermanifolds. Odd parameters}\label{OddComp}

Objects in the category of smooth supermanifolds are in one-to-one correspondence with vector bundles; this is proved in \cite{Ga, Bat}, \cite[Subsection~ 4.1.3]{MaG}. P.~Green, see \cite{G}, was the first to observe that in the category of complex-analytic supermanifolds there are more objects than there are vector bundles;  see also \cite{Pal} and Palamodov's chapter in \cite[Ch.4, \S4, Sections 6--9]{Ber1} reproduced in \cite[Ch.3, Th.2, p.126]{Ber2}. These Palamodov's chapters expound a~short note \cite{Pal} submitted for publication 2 months after~ \cite{G} was submitted, but published much later than \cite{G}. This topic attracted a~new interest with the discovery by Donagi and Witten of non-splitness of the moduli spaces of super Riemann surfaces, see \cite{DW}. For an overview of results on non-split supermanifolds obtained by Onishchik and his students, see \cite{Lob}; see also \cite[Ch.4,\S2, Prop.~8 (a), p. 190]{MaG} and \cite[Example 3.3.1 (1)]{Va2}. All papers on non-split supermanifolds published so far ignored odd parameters of deformations. Let me consider the odd parameters.

\subsection{The (almost) complex supermanifolds}\label{ssComplex} Let the ground field be $\Cee$, and $V$ an $m$-dimensional space. The sheaf
$ \mathcal F_{n|m}$
of germs of $\Lambda^{\bcdot}(V)$-valued holomorphic functions on
$\Cee^n$ defines the supermanifold~$\cC^{n|m}$. We consider the elements of $\mathcal F_{n|0}$  
even,  the elements $\xi_j$ of a~ basis of $V$  odd. A \textit{superdomain} in $\cC^{n|m}$ is a~ ringed space of the form $(U,\mathcal F_{n|m}|_U)$,
where $U$ is an open domain in~$\Cee^n$.
 An almost \textit{complex-analytic supermanifold} of dimension $n|m$ is a~
ringed space locally isomorphic to a
superdomain in~$\cC^{n|m}$. Thus, if $\cM=(M, \cO_\cM)$ is a~supermanifold, then for
any point $p_0\in M$ there exist a~neighborhood $U$ of $p_0$ in $M$
and an isomorphism of the ringed space $(U,\mathcal O|_U)$ onto a~superdomain
$(\widetilde U,\mathcal F_{n|m}|_{\widetilde U})$,  
where $\widetilde U$ is an open domain in $\Cee^n$, called a~\textit{chart}  in~$\cC^{n|m}$. Let $x_1,\dots ,x_n$ be coordinates in~$\Cee^n$.
Identifying $(U,\mathcal O|_U)$ with the superdomain by means of a~chart,
we get the sections $x_i$ for $i = 1,\dots ,n$, and $\xi_j$ for $j = 1,\dots ,m$ of
$\mathcal O|_U$ called  \textit{local coordinates} on $(U,\mathcal O|_U)$. 

For \textit{complex supermanifolds} and recently distinguished notion of \textit{real-complex supermanifolds}, such as Minkowski supermanifolds, see \cite{BGLS}.

\subsection{The retract of the supermanifold}\label{ssRetract}
Let $(M,\mathcal O_\cM)$ be a~supermanifold, $\mathcal I\subset\mathcal O_\cM$ a~subsheaf of
ideals; let $\mathcal O_\cM\tto \mathcal O_\cM/\mathcal I$ be the natural projection. Setting
\[
N = \{x\in M \mid \varphi(x) = 0\ \ \text{~~for all~~}\varphi \in\mathcal I\},\ \
\mathcal O_\cN = (\mathcal O_\cM/\mathcal I)|_N,
\]
we get the ringed space $(N,\mathcal O_\cN)$  called a~\textit{subsupermanifold} of $(M,\mathcal O_\cM)$. If
the sheaf $\mathcal I$ is generated, over an open set $U\subset M$, by
its
homogeneous sections $\varphi _1,\ldots,\varphi _s$, then it is usual to say that the
subsupermanifold is determined in $U$ by the system of equations $\varphi _i = 0$, where $i = 1,\ldots,
s$.

To any supermanifold
$\cM=(M, {\mathcal O}_\cM)$ the following split supermanifold corresponds. Consider the filtration
\begin{equation}
{\mathcal O}_\cM :=  {\mathcal J}^{0}\supset{\mathcal J}^{1} \supset {\mathcal J}^{2}
\supset \ldots\label{(1.12)}
\end{equation}
of $\mathcal O_\cM$ by the powers of the subsheaf of ideals ${\mathcal J}={\mathcal J}^{1}$ generated by the odd elements. The associated graded sheaf
$$
\gr{\mathcal O}_\cM=\bigoplus_{k\ge 0}\gr^k\mathcal O_\cM,
$$
where $\gr^{k}{\mathcal O}_\cM = {\mathcal J}^{k}/{\mathcal J}^{k+1}$,
gives rise to the split supermanifold $(M,\gr{\mathcal O}_\cM)$ 
called the \textit{retract} of the supermanifold~$\cM$. 
           
The quotient modulo the subsheaf of ideals $(\mathcal O_\cM)_\text{nil}$ generated by all nilpotents of $\mathcal O_\cM$ determines the \textit{reduction} whose result is conventionally denoted
$\cM_{\text{red}}$. The quotient modulo the subsheaf $\mathcal J^1$  generated by \textbf{only odd} nilpotents determines the \textit{odd reduction},  whose result Manin denoted 
$\cM_{\rd}$, of $(M,\mathcal O_\cM)$; thus, $\cM_{\rd}$ is  a~submanifold of $(M,\mathcal O_\cM)$.

Thus, we see that with any (almost) complex supermanifold $\cM$ two objects of the classical
complex analytic geometry are associated: the (almost) complex manifold $(M,\mathcal O_M)$
and the holomorphic vector bundle~$\mathbf E$ over $(M,\mathcal O_M)$ with fiber $V$ corresponding
to the sheaf $\mathcal E={\mathcal J}^{1}/{\mathcal J}^{2}$. It turns out that the complex supermanifold $\cM$ is not, in general, determined
by these two objects up to an isomorphism: there exist non-split
supermanifolds. 

The obstructions to splitness were considered and calculated, so far, for transition functions of the form \eqref{trFunEv}; e.g., see \cite{Pal, G}, \cite[Ch.4, \S4, Sections 6--9]{Ber1}, \cite[Ch.3, Th.2, p.126]{Ber2},  \cite[Ch.4,\S2, Prop.~8 (a), p. 190]{MaG},  \cite[Example 3.3.1 (1)]{Va2}, \cite{BO, Bash}. In this --- \textbf{restricted} --- approach (considering only transition functions of the form \eqref{trFunEv}), the first obstructions to splitness are given by the terms of degree $>1$ in the $\xi_i$ in the expansion of $\varphi^*(x_i)$ and $\varphi^*(\xi_j)$.  If the $x_i$ and $\xi_j$ are arbitrary local coordinates of $\cM$ in a~neighborhood $U
\subset M$, then $X_i = x_i +\mathcal J^2$ and $ \Xi_j = \xi_j +\mathcal J^3$ are splitting
local coordinates of $(M,\gr\mathcal O_\cM)$ in $U$, and one gets the transition
functions between these splitting coordinates, if one takes the terms of
degree 0 (respectively 1) in $\xi_j$ in the transition functions $\varphi _i$
(respectively $\psi_j$) for $\cM$, see eq.~\eqref{trFunEv}. 

In the restricted approach, the \textbf{first} obstructions to splitness are represented by cocycles in the \v{C}ech cohomology $H^1(M; \cT\otimes \cE^2(\mathbf E))$, where $\cE^i(\mathbf E)$ is the sheaf corresponding to the $i$th exterior power $E^i(\mathbf E)$ of the bundle $\mathbf E$, and $\cT$ is the tangent bundle, see \cite[Ch.4, \S2, Prop.~9, p. 191]{MaG} and Palamodov's chapter in \cite[Ch.4, \S4, Sections 6--9]{Ber1} reproduced in \cite[Ch.3, Th.2, p.126]{Ber2}. 

The \textbf{higher} obstructions to splitness in the restricted approach are represented by cocycles in the \v{C}ech cohomology $H^1(M; \cT \otimes \cE^{2i}(\mathbf E))$ for $i>1$, see  \cite[Ch.3, \S7, Th.2, p.135]{Ber2}; these obstructions describe even obstructions to splitness in the non-restricted approach as well, but miss odd obstructions.

In particular, this restricted approach (ignoring odd parameters, as in eq.~\eqref{trFunEv}) immediately leads one to the conclusion that ``any complex supermanifold of superdimension $m|1$ is split'' since $E^2(\mathbf E)=0$, see \cite[Ch.4,\S2, Prop.~8 (a), p. 190]{MaG} and \cite[Example 3.3.1 (1)]{Va2}. In \cite{MaG,Va1,Va2}, various deformation problems were considered over a~ supervariety of parameters, but --- inexplicably --- not in such generality in this particular case: in the study of non-splitness.

To incorporate odd parameters in a simplest case, recall the following well-known facts, see, e.g., \cite[\S~1.1]{OSS}. Let us cover ${\mathbb{CP}}^1$ by two affine charts $U_0$ and $U_1$ with local coordinates
$x$ and $y=\frac{1}{x}$, respectively. Let $\xi$ and $\eta$ be basis sections of  line bundle ${\mathbf L}_{k}$ over~
$U_0$ and $U_1$, respectively, such that
%the transition functions 
in $U_0\cap U_1$ we have (up to a non-zero constant factor)
\[ \begin{array}{l}
%y  ={x}^{-1},\ \ \
\eta  ={x}^{k}\xi. 
\end{array} 
\]
The bundle ${\mathbf L}_{k}$ is said to be \textit{of degree}~$k\in\Zee$. The sheaf of sections of this bundle is designated~$\cO(k)$. In particular, since $dy=-x^{-2}dx$, then the cotangent sheaf of 1-forms $\Omega^1$ is $\cO(-2)$, and hence %($\partial_y=x^{2}\partial_x$) 
the dual to~$\Omega^1$ tangent sheaf $\cT$ is $\cO(2)$; clearly, $\cO:=\cO(0)$ is the structure sheaf (of functions). We also have the following formulas, see \lq\lq Bott's formulas'' in \cite[p.4]{OSS}:
\be\label{facts}
\begin{array}{l}
\cO(a)\otimes\cO(b)\simeq \cO(a+b);\\

\dim H^0(\Cee\Pee^{1}; \cO(a))=\begin{cases}a+1&\text{if $a\geq0$}\\
0&\text{otherwise};\end{cases}\\
\text{$\xi, x\xi,\dots, x^{a}\xi$ is a~basis of the space of sections over $U_0$};\\

\end{array}
\ee
\be\label{facts1}
\begin{array}{l}
\dim H^1(\Cee\Pee^{1}; \cO(a))=\begin{cases}|a|-1&\text{if $a\leq-2$}\\
0&\text{otherwise};\end{cases}\\\text{for a~basis of $H^1(\Cee\Pee^{1}; \cO(a))$ we can take the classes of the}\\
\text{following cocycles on $U_0\cap U_1$:}\\
\begin{cases}x^{-1}\xi, x^{-2}\xi,\dots, x^{1-|a|}\xi&\text{if $a< -2$}\\
x^{-1}\xi&\text{if $a=-2$}.\end{cases}\\
\end{array}
\ee

Now, if we consider transition functions taking odd parameters into account, then eq.~\eqref{trFun} shows that there is a~linear in $\xi$ term in the new coordinates~$\varphi(x_i)$. The \textbf{odd} obstructions to splitness are represented by cocycles in $H^1(M; \cT\otimes \cE^{2i+1}(\mathbf E))$ for~$i\geq 0$.

In the case of the linear bundle $\mathbf E={\mathbf L}_{k}$ over $M=\Cee\Pee^{1}$, the corresponding obstructions belong to (recall that $E^1(V):=\Pi(V)$)
\be\label{*}
H^1(M; \cT\otimes \cE^1(\mathbf E))=H^1(M; \cO(2)\otimes \Pi(\cO(k)))\simeq H^1(M; \Pi(\cO(k+2))).
\ee

Eqs.~ \eqref{facts}, \eqref{facts1} and \eqref{*} imply the following theorem, where $\GL(V)=\GL(1)$. 

\sssbegin[(Odd parameters of the non-split $1|1$-dimensional superstring whose underlying  space is $\Cee\Pee^{1}$)]{Theorem}[Odd parameters of the non-split $1|1$-dimensional superstring whose underlying space is $\Cee\Pee^{1}$]\label{5.2.1} The non-split supermanifolds $\cM$ of complex superdimension $1|1$, a.k.a. superstrings, with underlying manifold $\Cee\Pee^{1}$, whose retract corresponds to the line bundle $\mathbf E={\mathbf L}_{k}$ with fiber $V$, \textbf{considered over a~base $\Spec C$ with $C_\od\neq 0$}, are described by $\GL(V)$-orbit in the space $H^1(M; \Pi(\cO(k+2)))$ --- its projectivization \mathstrike{, \text{the number of orbits being equal to}} \mathstrike{\text{$|k+2|-1$}}
for $k\leq -4$, whereas $\cM$ is split for~$k>-4$. \end{Theorem}

\subsection{On the claims of Bettadapura, see \cite{BK1, BK2, BK3}}\label{BK} Although Bettadapura sometimes uses standard notation in nonstandard ways (e.g., his $\Pi E$ is an \textit{ad hoc} notation, which in the context of super geometry causes confusion), writes that something ``will'' whereas it already ``is'', Palamodov's technique published in the chapter written by Palamodov (\cite[Ch.4, \S4]{Ber1}) is misattributed to Berezin, some symbols are left undefined (e.g., $E^{\text{v}}$), etc., etc., and claims are formulated assuming the change of coordinates is performed via formulas of the form \eqref{trFunEv}, not \eqref{trFun}, it is possible to decipher that 
some statements, correct over $\Spec\Cee$, become incomplete if we take into account odd parameters (e.g., \cite[Lemma~2.4]{BK1} --- a~copy of the corresponding claims in \cite{MaG} and \cite{Va2}, although these references are not given). These claims are at variance with Theorem~\ref{5.2.1}, the other statements \textit{do} take the odd parameters into account and contain ostensibly new examples of non-split super Riemann surfaces  $\cS$, but in reality they do not describe odd parameters of deformations of $\cS$, but rather an incomplete study of deformations of $\cS\times \cC^{0|n}$.
%\AO{I~softened \lq\lq wrong'' to \lq\lq correct over $\Spec\Cee$, but incomplete once odd parameters are taken into account'', for the same reason as in Theorem~\ref{5.2.1}: as stated over $\Spec\Cee$ those claims are right, and it is only for $C_\od\neq0$ that they fail to tell the whole story. Also deleted a~stray closing parenthesis after \lq\lq Theorem~\ref{5.2.1}''.}

Observe that, like many, but fortunately not everybody, Bettadapura applies the term ``super Riemann surface'' to $1|1$-dimensional supermanifolds with a~ contact structure, whereas it seems natural to consider \textbf{several} types of ``super Riemann surfaces'' of dimension $1|N$ either with or without contact structure but, perhaps, preserving a volume element, as I intend to show in a paper in preparation in which the upper bounds on $N>1$ are justified (explained) from a~certain point of view.

\section{Finite-dimensional simple Lie superalgebras considered naively}\label{LieListNaive} %, and their relatives}

Certain simple Lie superalgebras  (or their relatives, like $\fgl(n|n)$, $\fpgl(n|n)$, and $\fsl(n|n)$)  were easy to describe meaningfully, by some properties of these algebras. Even before a~matrix realization of $\fosp(m|2n)$ and its interpretation as a~Lie superalgebra preserving a~non-degenerate even bilinear form were discovered it was natural to call it ``ortho-symplectic'', just by looking at its even part. Two series of simple Lie superalgebras were considered ``strange'' until the interpretations of their non-simple relatives as preserving a~tensor was discovered, hence their first notation and certain awkward names (e.g., \lq\lq di-spin algebra of Mitchell---Gell-Mann---Radicati\rq\rq) were soon abandoned; certain confusing names --- $F(4)$ and $G(3)$ instead of \textit{ad hoc} but not confusing (with the Lie algebra $\ff_4$ and some object on the line of exceptional Lie algebras on the Vogel plane, see \cite{MSV}, respectively) names given by their discoverer, Kaplansky) still remain in some papers. The notation ``$\fpsl$'' was also being avoided for a~long time, although it was already usual in the modular setting.

\subsection{Linear (matrix) Lie superalgebras} Let $\Size:=(p_1, \dots, p_{|\Size|})$, where $|\Size|:=\dim V$, be an ordered
collection of parities of the basis vectors of a superspace $V$ for which we take only vectors \textit{homogeneous with respect to parity}. The \textit{general linear} Lie
superalgebra of all supermatrices of size $\Size$ corresponding to linear operators in the superspace $V=V_{\bar 0}\oplus V_{\bar 1}$ over the ground field $\Kee$ is denoted by~$\fgl(\Size)$. 
Usually, for the \textit{standard} (a~simplest from a~certain point of view) format, $\fgl(\ev,
\dots, \ev, \od, \dots, \od)$ is abbreviated to $\fgl(\dim V_{\bar
0}|\dim V_{\bar 1})$. Any non-zero supermatrix in $\fgl(\Size)$ can
be uniquely expressed as the sum of its even and odd parts; in the
standard format, this is the following block expression; on non-zero summands the parity is defined:
\[
\mmat A,B,C,D,=\mmat A,0,0,D,+\mmat 0,B,C,0,,\quad
 p\left(\mmat A,0,0,D,\right)=\ev, \; p\left(\mmat 0,B,C,0,\right)=\od.
\]
The \textit{trace} on the Lie superalgebra $\fg$ (or \textit{supertrace} for emphasis) is any map ${\tr:\fg\tto\Kee}$ that vanishes on the first derived Lie (super)algebra~$\fg'$. In particular, on $\fgl (\Size)$, there is just one (up to a~non-zero scalar) such map; it is given by the formula
\be\label{str}
\str (X):=\sum (-1)^{p_{i}}X_{ii}.
\ee
Thus, in the standard format, 
\[
\str \begin{pmatrix}A&B\\ C&D\end{pmatrix}=\tr A- \tr D.
\]

Observe that for the Lie superalgebra $\fgl_\cC(a|b)$ over a~supercommutative superalgebra $\cC$, i.e., for supermatrices with homogeneous (with respect to parity) elements in $\cC$, we have
\be\label{str1}
\begin{array}{l}
\str X=\tr A- (-1)^{p(X)}\tr D\text{~~for any $X=\begin{pmatrix}A&B\\ C&D\end{pmatrix}$,}\\
\text{where $p(X)=p(A_{ij})= p(D_{kl})=p(B_{il})+\od=p(C_{kj})+\od$},\\
\end{array}
\ee
so on odd supermatrices with entries in $\cC$ such that $\cC_\od\neq 0$, the expression of supertrace coincides with the conventional trace of a matrix.

Since $\str [x, y]=0$, the supertraceless
supermatrices span a~Lie subsuperalgebra called \textit{special linear} and denoted~$\fsl(\Size)$.

There are, however, at least two super versions of $\fgl(n)$, not
one; for reasons, see \cite[Ch1, Ch.7]{LSoS}, where the super version of the Schur lemma and classification of central simple (finite-dimensional) associative superalgebras are considered. The natural (from the Schur lemma's point of view) other version --- $\fq(n)$ --- is called the \textit{queer}
Lie superalgebra and is defined as the one that preserves --- over an algebraically closed field $\Kee$ of characteristic $p\neq 2$ --- the
complex structure given by an \textit{odd} operator~$J$, i.e.,
$\fq(n)$ is the supercentralizer $\mathrm{C}(J)$ of~$J$:
\[
\fq(n):=\mathrm{C}(J)=\{X\in\fgl(n|n)\mid [X, J]=0 \}, \text{ where }
J^2=-\id.
\]
It is clear that over an algebraically closed field,  we can by a~change of basis reduce $J$ to the normal form (shape)
$J_{2n}$ in the standard
format,  see eq.~\eqref{q},
and then the elements of $\fq(n)$ take the form
\begin{equation}\label{q}
\fq(n):=\left \{(A,B):=\mat {A&B\\B&A}, \text{~~where $A, B\in\fgl(n)$ and $J_{2n}:=\mat {0&1_n\\-1_n&0}$}\right\}.
\end{equation}
(Over any 
field $\Kee$, instead of $J$ we can take any odd operator $K$ such that ${K^2=a\id_{n|n}}$, the multiple of the identity operator, where $a\in \Kee^\times$. If $\Kee$ is algebraically closed, the Lie superalgebras $\mathrm{C}(K)$ --- supercentralizers of $K$ --- are isomorphic for distinct $a$; if $p=2$, it is natural to select~$K^2=\id$.)

On $\fq(n)$, the odd trace, nowadays called \textit{queertrace}, is defined: 
\[
\qtr\colon (A,B)\longmapsto \tr B. 
\]
Denote by $\fsq(n)$ the Lie superalgebra
of \textit{queertraceless} matrices, first described in \cite{BL2} together with the odd (queer) determinant of the $C$-points of the corresponding supergroup, see a~more accessible \cite[Ch.1]{LSoS}.

Clearly, $\fgl$ and $\fq$ correspond to the two cases of the super version of Schur's lemma over an algebraically closed field: \textbf{the $\fgl$ case}: an irreducible module $M$ over a~collection $S$ of homogeneous operators can be \textit{absolutely irreducible}, i.e., have no proper invariant subspaces, then the only operator commuting with $S$ is a~scalar, and \textbf{the $\fq$ case}:  $M$  can have an invariant subspace, which is not a~sub\textbf{super}space, then the superdimension of the module is of the form $n|n$ and an odd operator $K$ interchanges the homogeneous components of the module.

\subsubsection{Supermatrices of operators}
To a~ linear map of superspaces ${F: V\tto W}$ there corresponds the
dual map $F^*:W^*\tto V^*$ between the dual superspaces. In bases
consisting of homogeneous vectors $v_{i}\in V$ of parity $p(v_i)$, and $w_{j}\in W$ of parity $p(w_j)$, the formula
$F(v_{j})=\mathop{\sum}_{i}w_{i}X_{ij}$ assigns to $F$ the
supermatrix~$X$. In the dual bases, the \textit{supertransposed}
matrix $X^{st}$ corresponds to $F^*$:
\[
(X^{st})_{ij}:=(-1)^{(p(v_{i})+p(w_{j}))p(w_{j})}X_{ji}.
\]
In the standard supermatrix format, we have
\[
\begin{footnotesize}
X=\begin{pmatrix}A&B\\C&D\end{pmatrix}\longmapsto X^{st}:=
\begin{pmatrix}A^t&(-1)^{p(X)}C^t\\-(-1)^{p(X)}B^t&D^t\end{pmatrix}=\begin{cases}\begin{pmatrix}A^t&C^t\\-B^t&D^t\end{pmatrix}&\text{if $p(X)=\ev$},\\
\begin{pmatrix}A^t&-C^t\\B^t&D^t\end{pmatrix}&\text{if $p(X)=\od$}.\end{cases}
\end{footnotesize}
\]

\subsubsection{Supermatrices of bilinear forms}\label{sssBilMatr}
 Having selected a~basis (by definition consisting of vectors homogeneous with respect to parity) of the superspace $V$, we define the Gram matrix $B=(B_{ij})$ of the
bilinear form $B^f$ on $V$ by the formula
\be\label{martBil}
B_{ij}:=(-1)^{p(B)p(v_i)}B^f(v_{i}, v_{j})\text{~~ for the basis vectors $v_{i}\in V$.}
\ee
This formula for the Gram matrix of $B^f$ allows us
to identify any~bilinear form $B\in B(V, W)$, where $V$ and $W$ are superspaces, with an operator, an element of $\Hom(V, W^*)$, see \cite[Ch.1]{LSoS}.

Recall that the \textit{upsetting} of bilinear forms
$u\colon\Bil (V, W)\tto\Bil(W, V)$ is given by the formula
\be\label{susyB}
u(B^f)(w,v):=(-1)^{p(v)p(w)}B^f(v, w)\text{~~ for any $v \in V$ and $w\in W$.}
\ee

Let now $W=V$, and $\Bil (V):=\Bil (V, V)$. The shape of the Gram matrix $B^u$ of a~homogeneous form $u(B^f)$ in the standard format of $V$ is as follows
\be\label{BilSy}
B^{u}:=
\mmat R^{t},(-1)^{p(B)}T^{t},(-1)^{p(B)}S^{t},-U^{t},\text{~~ for $B=\mmat R,S,T,U,$,}
\ee

The form
$B^f$ is said to be \textit{symmetric} if $B^{u}=B$, and \textit{anti-symmetric} if~$B^{u}=-B$. 
In particular, the form on $\fgl(V)$ (resp. $\fq(V)$) given by the respective trace is symmetric:
\[
\begin{array}{l}
(X,Y):=\str(XY) \text{~~for any $X,Y\in \fgl(V)$}\\ (\text{resp. $(X,Y):=\qtr(XY)$ for any $X,Y\in \fq(V)$)}.
\end{array}
\]

Clearly, the \textit{upsetting} $B\mapsto B^{u}$ of Gram matrices of bilinear forms is \textit{not} supertransposition.

Observe
that \textbf{the passage from $V$ to $\Pi (V)$ turns every symmetric
form $\cB$ on $V$ into an anti-symmetric one $\cB^\Pi$ on $\Pi (V)$  and every anti-symmetric $\cB$ into a~symmetric $\cB^\Pi$  by setting}
\[
\text{$\cB^\Pi(\Pi(x), \Pi(y)):=(-1)^{p(B)+p(x)+p(x)p(y)}\cB(x,y)$ for any $x,y\in V$}.
\]

Observe that there are no ``supersymmetric'' bilinear forms; the property 
\be\label{braid}
B^f(w,v):=(-1)^{p(v)p(w)}B^f(v, w)\text{~~ for any $v, w \in V$}
\ee
which is sometimes given as their definition just reflects braiding in the category of superspaces and has nothing to do with the symmetry of bilinear forms. (Therefore, the title of  \cite{BKLS} should be corrected. The same confusing term is used in \cite{BS}, where inhomogeneous with respect to parity bilinear forms are classified on superspaces of dimension $\leq 7$ in an attempt to define a~generalized mixture of the Heisenberg superalgebra and what in  \cite{BGLLS} was called the \textit{antibracket superalgebra}. The mainstream supersymmetry theory ignores the inhomogeneous bilinear forms since they are not preserved by Lie superalgebras; but it is interesting, nevertheless, to explore their hidden depths, cf. Subsection~\ref{Uma}; e.g., what is the ``something'' that preserves them? What linear transformations of the superspace $V$ (only even or any) is reasonable in some (which?) sense to consider in the classification of inhomogeneous bilinear forms on~$V$?) 

Most popular normal shapes of the (Gram matrices of the) even non-degenerate symmetric
form are the ones which in the standard format are as
follows:
\[
\begin{array}{l}
B_{ev}(m|2n)= \diag(1_m, J_{2n}):=\mmat 1_m,0,0,J_{2n}, \text{~or~} \diag(A_m,J_{2n}):=\mmat A_m,0,0,J_{2n},,\\ \text{where
$J_{2n}= \antidiag(1_n,-1_n):=\mmat 0, 1_n,-1_{n},0,$ and $A_m=\antidiag(1,\dots, 1)$}.\\
\end{array}
\]

The Lie superalgebra $\faut (B)\subset \fgl(\Size)$ that preserves the Gram matrix $B$ of the form ${B^f\in\Bil(V)}$ consists of the supermatrices $X\in\fgl(\Size)$ such that
\[
X^{st}B+(-1)^{p(X)p(B)}BX=0\quad \text{for an homogeneous matrix
$B\in\fgl(\Size)$}.
\]

The usual notation for $\faut (B_{ev}(m|2n))$ is $\fosp(m|2n)$;
sometimes one writes more explicitly, $\fosp^{sy}(m|2n)$. Observe
that the anti-symmetric non-degenerate \textbf{even} bilinear form is preserved by the
``symplectic-orthogonal'' Lie superalgebra, $\fsp\fo(2n|m)$ or $\fosp^{a}(m|2n)$, which is
isomorphic to $\fosp^{sy}(m|2n)$, but has a~different matrix realization.

A~non-degenerate \textbf{symmetric} odd bilinear form $B_{odd}(n|n)$ can
be reduced to a~normal shape whose matrix in the standard format is
$J_{2n}$, see \eqref{BilSy}, \textbf{not} $\Pi_{2n}:=\antidiag(1_n,1_n)$ which is \textbf{anti}symmetric, see \cite{LSoS}, contrary to a~hasty conclusion induced by the ``even'' experience. The usual notation for
$\faut (B_{odd}(\Size))$ is~$\fpe(\Size)$.
The passage from $V$ to $\Pi (V)$ establishes an isomorphism
$\fpe^{sy}(\Size)\cong\fpe^{a}(\Size)$. These isomorphic Lie
superalgebras are nowadays called, as A.~Weil suggested, \textit{periplectic}.

The traceless (special) subalgebra of  $\fpe(n)$ is denoted $\fspe(n)$; it is simple for~$n>2$.

For a~large class of Lie superalgebras either simple, or relatives of simple, Kac  introduced the notion of a~Cartan matrix, see \cite{K2} with improvements of the definition in \cite{CCLL, BGL, BLS, BLLoS}. If $p\neq 2$, then neither periplectic superalgebras nor their simple relatives have Cartan matrices; this is not so for $p=2$, see the paper~\cite{BGL} containing the~classification of all finite-dimensional modular Lie superalgebras with indecomposable Cartan matrix over any algebraically closed field of positive characteristic.

If $\fc$ is a~Lie algebra of scalar matrices, and $\fg\subset \fgl (n|n)$ is a~Lie subsuperalgebra containing $\fc$, then the \textit{projective}
Lie superalgebra of type $\fg$ is $\fpg:= \fg/\fc$, e.g., $\fpsl(n|n)$ and~$\fpsq(n)$. 

\subsection{Exceptional finite-dimensional simple Lie superalgebras} There are 3 of them: $\fosp_a(4|2)$, $\fag(2)$ and~$\fab(3)$. I denote them following meaningful Serganova's notation from the paper \cite{Se4}, where she proved that the notion of roots generalized to non-degenerate \textbf{pseudo}-euclidean space instead of euclidean, leads --- unexpectedly --- to Lie superalgebras $\fg(A)$ with indecomposable matrix~$A$. A~tempting \textbf{Open problem}: superise the result of \cite{P} by considering non-degenerate but not sign-definite sesquilinear form --- an analog of Serganova's result for complex reflections. see \cite{Se4}.

Kaplansky was the first to discover all three exceptional simple finite-dimensi\-o\-nal Lie superalgebras over $\Cee$ (\cite{Kapp}). They are of the form $\fg(A)$, i.e., have a~ Cartan matrix $A$, and thus are described in \cite{BGL} over fields of any characteristic $p$ (non-existing for some values of $p$), where the description as pairs $(\fg_\ev, \fg_\od)$, although snubbed at in Subsection~\ref{clumsy} above as ``clumsy'', is also given because it is useful in certain calculations, see also \cite{CCLL}. For \textit{interpretations} of the exceptional simple Lie algebras, see  \cite{Eld, Eldi} and \cite{Sud}, see also Subsection~\ref{ssExc}.

\subsection{Finite-dimensional simple vectorial Lie superalgebras} Define partial derivatives in $\Lambda^{\bcdot}(n)$ by setting $\partial_{\xi_i}(\xi_j):=\delta_{ij}$ and considering only left action of the derivations on functions. Let $\partial_i:=\partial_{\xi_i}$. Then, setting $\deg \xi_i=1$ for all $i$, we get the \textit{standard} $\Zee$-grading of 

$\bullet$ the \textit{general vectorial Lie superalgebra} 
\[
\fvect(0|n):=\left\{\sum f_i\del_i\mid f_i\in\Lambda^{\bcdot}(n)\text{~~for all $i$}\right\}=\mathop{\oplus}\limits_{-1\leq i\leq n-1}\fg_i.
\]

$\bullet$ The \textit{divergence} (with respect to the volume element $\vvol(\xi)$ in ``coordinates'' $\xi$ which means that the Lie derivative $L_D$ along the vector field $D$ multiplies $\vvol(\xi)$ by $\Div(D)$) defined by the formula (clearly, $p(\del_i)=\od$ for all $i$)
\[
\Div \left(\sum f_i\del_i\right):= \sum (-1)^{p(f_i)}\del_if_i.
\]

The \textit{special} or \textit{divergence--free} vectorial Lie superalgebra is
\[
\fsvect(0|n):=\{D\in \fvect(0|n)\mid \Div (D)=0\}=\{D\in \fvect(0|n)\mid L_D(\vvol(\xi))=0\}.
\]

Kac was the first to observe (in analogy with characteristic $>0$, cf. \cite{W}) that for $n=2k$ any \textit{volume form} $f\vvol(\xi)$, where $f\in\Lambda^{\bcdot}(\xi)$, can be reduced by a~change of indeterminates (an automorphism of $\Cee[\xi]$) to one of the two \textit{normal shapes}: with either $f=1$ or with $f=(1+\prod \xi_i)$, see \cite{K1.5}, \cite[Section~3.3.1]{K2}, where the volume forms (a particular case of \textit{integrable forms}, see \cite{BL}) are erroneously considered as a~particular case of differential forms. Kac formulated this even before the notion of the volume form was correctly defined. The same argument works for functions in $n=2k+1$ indeterminates with values in $\Cee[\tau]$, where $p(\tau)=\od$, the normal shapes being with either $f=1$ or with $f=(1+\tau\prod \xi_i)$. Accordingly, the \textit{deformed special} or \textit{divergence--free} vectorial Lie superalgebra is
\[
\widetilde{\fsvect}(0|2k):=\{D\in \fvect(0|2k)\mid L_D((1+t\prod \xi_i)\vvol(\xi))=0\}.
\]
Any two such Lie superalgebras for $t, t'\in\Cee^\times$ are clearly isomorphic, so we assume~$t=1$.

Considering the \textit{supergroup of automorphisms} $\Cee[\xi_1,\dots, \xi_n]$ in terms of the functor 
\[
\cC\longmapsto \Aut^\ev_\cC(\cC\otimes\Cee[\xi]),
\]
where $\cC$ is a~supercommutative superalgebra (cf. with the description of a~ group scheme in \cite{MaAG}), we arrive at the following normal shapes of the volume forms $f\vvol(\xi)$:
\[
f=\begin{cases}1;\\
(1+t\prod \xi_i), \text{~~where $t\in\cC_\ev^\times$}&\text{for $n$ even},\\
(1+\tau\prod \xi_i), \text{~~where $\tau\in\cC_\od$}&\text{for $n$ odd},
\end{cases}
\]
and hence to an odd parameter $\tau$ and the \textit{deformed special} or \textit{divergence--free} vectorial Lie superalgebra 
\[
\widetilde{\fsvect}(0|2k+1):=\{D\in \Cee[\tau]\otimes \fvect(0|2k+1)\mid L_D((1+\tau\prod \xi_i)\vvol(\xi))=0\}.
\]
Any two such Lie superalgebras for odd parameters $\tau,\tau'$ are, clearly, isomorphic.

$\bullet$ 
Define the space of \textit{differential} (a.k.a. \textit{exterior}) forms $\Cee[\xi, d\xi]$ by setting $p(d\xi_i)=\ev$ for all $i$, and letting the $d\xi_i$ commute with the $\xi_j$ and $d\xi_j$; define the \textit{exterior differential} by the formula (together with the Leibniz and Sign Rules)
\[
d(f):=\sum d\xi_i\frac{\partial f}{\partial \xi_i}.
\] 

Define the action of \textit{Lie derivative} along the vector field $D$ on the space of differential  forms $\Cee[\xi, d\xi]$ by setting
\[
L_D(df):=(-1)^{p(D)} d(D(f)).
\]

For a~normal shape of the \textit{symplectic}, i.e., a~non-degenerate and closed, 2-form  on $\cC^{0|n}$ we can take 
\[
\omega:=\sum (d\xi_i)^2.
\]
This was first claimed in \cite{Ln}, for the  first proof, see \cite{Sha}; for a~conceptual proof, see \cite{GLr, BGLS}; observe that in characteristic $p>0$ there are many normal shapes of the symplectic forms, see \cite{Sk} where this classification is performed  on manifolds; to superize Skryabin's result is an important \textbf{Open problem}. 

Let the \textit{Lie superalgebra of Hamiltonian vector fields} on $\cC^{0|n}$ be
\[
\fh(0|n):=\left\{D\in \fvect(0|n)\mid L_D(\omega)=0\right\}.
\]

Let the space of the \textit{Poisson} Lie superalgebra $\fpo(0|n)$ on $\cC^{0|n}$ be $\Lambda^{\bcdot}(n)$ with the \textit{Poisson bracket}
\[
\{f,g\}_{P.b.}:=(-1)^{p(f)}\sum \frac{\partial f}{\partial \xi_i} \frac{\partial g}{\partial \xi_i}.
\]
Let $H_-: f\longmapsto H_f$ be the projection $\fpo(0|n) \tto \fh(0|n)$ given by the formula
\[
H_f:=(-1)^{p(f)}\sum \frac{\partial f}{\partial \xi_i} \partial_{\xi_i}.
\]
Clearly, the following sequence is exact
\be\label{PoH}
0\tto \Cee\tto \fpo(0|n) \stackrel{H_{-}}{\tto} \fh(0|n) \tto 0.
\ee

\subsection{The list of finite-dimensional simple Lie
superalgebras considered naively} Recall that the Lie superalgebra
$\fg$ without proper ideals and of dimension~$>1$ is said to be
\textit{simple}. Eq.~\eqref{simple} lists all non-isomorphic (minding isomorphisms \eqref{isom}, \eqref{osp42symm}, \eqref{osp42symm1}) simple
finite-dimensional Lie superalgebras over $\Cee$; for first examples and
subcases of classification, see \cite{K0,K1, Kapp, K1C, SNR, Ho, Dj1, Dj2, Dj3, DjH}; for a
summary (ignoring odd parameters), see \cite{K2, Sch} with improvements cited in eq.~\eqref{simple}:
\begin{equation}\label{simple}
\begin{array}{l} \text{$\fsl(m|n)$ for $m> n\geq 1$, $\fpsl(n|n)$ for
$n>1$, $\fosp(m|2n)$ for $mn\neq 0$}\\
\text{$\fpsq(n)$ and $\fspe(n)$
for $n>2$;} \\
\text{$\fosp_\alpha(4|2)$ for $\alpha\neq 0, -1$ is given by any of its four Cartan matrices (\cite{BGL})} \\ \begin{footnotesize}
\begin{pmatrix}
0&1&-1-\alpha\\
-1&0&-\alpha\\
-1-\alpha&\alpha&0\\
\end{pmatrix} \ \begin{pmatrix}
2&-1&0\\
-1&0&\frac{1+\alpha}{\alpha}\\
0&-1&2\\
\end{pmatrix}\ \begin{pmatrix}
2&-1&0\\
-1&0&-\alpha\\
0&-1&2\\
\end{pmatrix}\ \begin{pmatrix}
2&-1&0\\
-1&0&\frac{1}{1+\alpha}\\
0&-1&2\\
\end{pmatrix}\end{footnotesize}
\\
\text{$\fa\fg(2)$ and
$\fa\fb(3)$ are also given by their Cartan matrices (see \cite{BGL});}\\
\text{for interpretations, see \cite{DWN,Sud} and very interesting papers \cite{Eld, Eldi};}\\
%\text{finite-dimensional vectorial Lie superalgebras, see Section~\ref{SS:2.5}};\\
\text{$\fvect(0|n)$ for $n>1$,} \\
\text{$\fsvect(0|n)$  for $n>2$ and
$\widetilde\fsvect(0|2n)$ for $n>1$,} \\
\text{$\fh'(0|n):=[\fh(0|n), \fh(0|n)]$ for $n>3$.}\\
\end{array}
\end{equation}

\textbf{Occasional isomorphisms}:
\be\label{isom}
\begin{array}{l}
\fsl(m|n)\simeq \fsl(n|m),\ \ \fvect(0|2)\simeq \fsl(1|2)\simeq\fosp(2|2),\\
\fspe(3)\simeq\fsvect(0|3),\ \ \fpsl(2|2)\simeq\fh'(0|4),\\
\fosp_1(4|2)\simeq\fosp_{-2}(4|2)\simeq\fosp_{-1/2}(4|2)\simeq\fosp(4|2).\\
\end{array}
\ee
The isomorphisms $\fosp_{\alpha}(4|2)\simeq \fosp_{\alpha'}(4|2)$ are generated by the
transformations (see \cite{BGL}):
\begin{equation}\label{osp42symm}
\alpha\longmapsto \alpha':= -1-\alpha\ , \qquad \alpha\longmapsto \alpha':=
\nfrac{1}{\alpha}\ ,
\end{equation}
so the other isomorphisms are 
\begin{equation}\label{osp42symm1}
\alpha\longmapsto \alpha':=
-\nfrac{1+\alpha}{\alpha}\ ,\quad\alpha\longmapsto \alpha':=
-\nfrac{1}{\alpha+1}\ ,\quad\alpha\longmapsto \alpha':= -\nfrac{\alpha}{\alpha+1}.
\end{equation}

In what follows, I'll list all deformation of simple finite-dimensional Lie superalgebras (see Theorems \ref{ThEv} and \ref{ThOd}), in particular, I'll prove the old conjecture: \textbf{The Lie superalgebras $\widetilde\fsvect(0|2n+1)$ for $n\geq 1$ are the only (bar $\fspe(3)$, see eqs.~\eqref{isom}) deformations of simple finite-dimensional Lie superalgebras with an odd parameter}.

\subsubsection{Exceptional simple Lie superalgebras: realizations}\label{ssExc} The irreducible modules of least dimension over $\fag(2)$, $\fab(3)$, and $\fosp_a(4|2)$ for $a$ generic are the adjoint ones. Therefore, for the best interpretation of these algebras, see \cite{Eldi}.

Consider $\fosp_a(4|2)$ for the exceptional values of~$a$. Kac was, probably, the first to establish their irreducibles of least dimension in an unpublished preprint \cite{Kpre}.
%\AO{The key of this \texttt{bibitem} was \texttt{K?} and printed as \lq\lq [K?]'', which looks like an unfinished reference; I~renamed it to \texttt{Kpre}. I~could not find this 1978 preprint in any database. If you have it on paper, a~scan on your page would let the reader check the claim on the least-dimensional irreducibles of $\fosp_a(4|2)$ for exceptional $a$.}

For $a=1$ (and hence for $a=-2$ and $a=-\frac12$), the Lie superalgebra $\fosp_a(4|2)$ is   isomorphic to $\fosp(4|2)$ and its irreducible module of least dimension is (up to the change of parity) the $(4|2)$-dimensional tautological module in which the non-degenerate even symmetric bilinear form is preserved.

For $\fosp_2(4|2)$ and $\fosp_3(4|2)$ taken in realization with respective Cartan matrices
\be\label{cm23}
\begin{pmatrix}
2&-1&0\\
-1&0&-2\\
0&-1&2\\
\end{pmatrix} \text{~~and~~} \begin{pmatrix}
2&-1&0\\
-1&0&-3\\
0&-1&2\\
\end{pmatrix}
\ee
--- the shape agreed upon above, normalized so that the isotropic row has $-1$ in the left
column, with $t=a=2$ and $3$ --- consider the weights of irreducibles over them relative
their Chevalley generators ${H_i:=[X_i^+,X_i^-]}$.

For $a=2$, as well as for $a=-3$,  $\frac12$, $-\frac32$, $-\frac13$ and $-\frac23$,  the irreducible module $V$ of least dimension is (up to the change of parity) the $\fosp_2(4|2)$-module with $\sdim V=(6|4)$ and an odd highest weight vector of weight $(1,-2,0)$ with respect to the Cartan subalgebra of~$\fg_\ev$.

For $a=3$, as well as for $a=-4$,  $\frac13$, $-\frac43$, $-\frac14$ and $-\frac34$,  the irreducible module $W$ of least dimension is (up to the change of parity) the $\fosp_3(4|2)$-module with $\sdim W=(6|8)$ and an even highest weight vector of weight $(2,-3,0)$ with respect to the Cartan subalgebra of~$\fg_\ev$.

Interestingly, the two exceptional algebras are somewhat similar in structure, see \cite{GL4}, being represented as a direct sum $\bigoplus $ --- as superspaces --- of a~subsuperalgebra $\fh$ and an $\fh$-module:
\[
\begin{array}{l}
\fag(2) = \fosp_3(4|2)\bigoplus  \Pi(W); \\ \fab(3) = (\fosp_2(4|2) \oplus  \fsl(2)) \bigoplus  (V \boxtimes  R(\omega_1)),  \end{array}
\]
where $R(\omega_1)$ is the tautological $\fsl(2)$-module and $A\boxtimes B$ is the modern notation of the tensor product of modules over different algebras, e.g., of the $\fg$-module $A$ by the $\fh$-module~$B$.

\section{Description of deformations}\label{defProof}

\ssbegin[(Classification of deformations with odd parameters)]{Theorem}[Classification of deformations with odd parameters]\label{ThOd} The only \textup{(observe $\fsvect(0|3)\simeq\fspe(3)$, see eq.~\eqref{isom})} simple finite-dimensional Lie superalgebras over $\Cee$ having deformations with odd parameters are $\fg=\fsvect(0|2n+1)$ for $n\geq 1$ and $\sdim H^2(\fg; \fg)=(0|1)$. The deformed superalgebra is $\widetilde{\fsvect}(0|2n+1)$.
\end{Theorem}

This theorem was conjectured in 1976 by Bernstein and me, and M.~Duflo as he told me in 1986; this conjecture was never published. In this section, I'll prove it. For the vectorial Lie superalgebras, I know only one method to prove this conjecture; it requires ingredients published in the paper \cite{GLP}. This section contains also the details of the proof and description of deformations with even parameters; almost all cases were partly described earlier by other researchers.

\ssbegin[(Classification of deformations with even parameters)]{Theorem}[Classification of deformations with even parameters]\label{ThEv} The only \textup{(minding isomorphisms \eqref{isom})} simple finite-dimensional Lie superalgebras over $\Cee$ having deformations with even parameters are $\fg=\fsvect(0|2n)$ for $n\geq 2$, $\fh'(0|m)$ for $m\geq 5$, and~$\fosp(4|2)$. In all these cases, $\dim H^2(\fg; \fg)=1$. The respective deforms are $\fg=\widetilde{\fsvect}(0|2n)$, $\fpsl(2^{k-1}|2^{k-1})$ for $m=2k$ and $\fpsq(2^{k-1})$ for $m=2k-1$, and $\fosp_\alpha(4|2)$. The deforms are isomorphic for any two non-zero values of parameter; for  isomorphisms of $\fosp_\alpha(4|2)$ for distinct $\alpha$, see eq.~\eqref{osp42symm}.
\end{Theorem}

\ssec{General facts on cohomology of Lie algebras and Lie superalgebras}
See \cite[pp.288--289]{Kn} and \cite[Ch.1]{Fu}. For correct answers of $H^{\bcdot}(\fg)$ for simple finite-dimensional Lie superalgebras $\fg$ over $\Cee$, see \cite{BoKN}, except for $\mathfrak{svect}(0|n)$, $\widetilde{\mathfrak{svect}}(0|n)$ and $\mathfrak{h}'(0|n)$ in which case to describe $H^{\bcdot}(\mathfrak{g})$ is an Open problem, cf. \cite{AF}. For the peculiarities of description of deformations in characteristic $p>0$, see \cite{BGL1, BLLS}; for general definitions, see a~ detailed discussion in \cite{BGL1}.

\ssec{Lie superalgebras $\fg(A)$ with indecomposable Cartan matrix $A$} Kac was the first to observe that some of simple Lie superalgebras have Cartan matrices, and that one Lie superalgebra can have several Cartan matrices, see \cite{K2}. For a~correct definition of Cartan matrices, and how to list all Cartan matrices of a~given Lie superalgebra, see \cite{CCLL, BLS, BLLoS}. 
For the Lie superalgebra $\fg=\fg(A)$, if $A$ is invertible \textbf{and the Killing form of $\fg$ does not vanish identically}, the proof of the fact ``${H^2(\fg;\fg)=0}$'' is the same as for simple Lie algebras (see \cite{Ru}) with the help of quadratic Casimir element corresponding to the even non-degenerate invariant symmetric bilinear form (even NIS) on $\fg(A)$; for the description of such forms on Lie superalgebras $\fg(A)$, see \cite{BKLS}.
%\AO{The condition in bold is new: without it the claim is false, see subsec.~\ref{KrugCas}, where the exceptions are listed and settled. Note that the invertibility of $A$ does \textbf{not} imply that the Killing form is non-zero: for a~Lie superalgebra, unlike for a~Lie algebra, the two conditions are unrelated.}

However, $\fpsl(n|n)$ has no Cartan matrix, while none of several Cartan matrices of $\fgl(n|n)$ is invertible.

If $A$ is not invertible, we have to consider the simple subquotient $\fg'(A)/\fc$ of $\fg(A)$, and the proof is still the same with the help of the quadratic Casimir element corresponding to the NIS induced from the NIS on~$\fg(A)$.

\sssec{B.~Kruglikov's corrections: the Casimir element dies exactly when the Killing form does}\label{KrugCas}
B.~Kruglikov observed (a~letter to me of 20 February 2025) that the argument
just given says nothing about $\fg=\fosp(2n+2|2n)$, where $n\geq 1$: for these
Lie superalgebras the Cartan matrix is invertible, whereas the Killing form is
identically zero. He is right, and the omission is not innocent: for
$\fosp(4|2)$, which is the case $n=1$, the conclusion is \textbf{false}, as
Theorem~\ref{ThEv} of this very paper says. Here is the precise state of affairs.

Recall that the \textit{Killing form} of $\fg$ is
${\kappa(x,y):=\str(\ad_x\ad_y)}$, and that on a~simple Lie superalgebra $\fg$ with an even NIS $B$
we have ${\kappa=cB}$ for a~constant $c$, since the space of even non-degenerate invariant
symmetric bilinear forms on $\fg$ is 1-dimensional by Schur's lemma (note that in certain situations, this is not so, see \cite{KLLS}).

\sssbegin{Lemma}[The Casimir element kills $\fg$ exactly when $\kappa=0$]\label{CasKil} Let $\fg$ be a~simple finite-dimensional Lie
superalgebra over $\Cee$ with an even NIS $B$, let $\Omega$ be the quadratic
Casimir element corresponding to $B$, and let ${\dim\fg_\ev\neq\dim\fg_\od}$. Let~${\kappa=cB}$. Then, $\Omega$ acts on the adjoint module by the scalar
\[
c\cdot\frac{\dim\fg_\ev+\dim\fg_\od}{\dim\fg_\ev-\dim\fg_\od}.
\]
In particular, $\Omega$ acts on $\fg$ by $0$ if and
only if~${\kappa=0}$.
\end{Lemma}

\begin{proof} As $\fg$ is simple, the $\fg$-module $\fg$ is irreducible, so
$\Omega$ acts on it by a~scalar $c'$, and hence
\[
\str_\fg\Omega=c'\str(\id_\fg)=c'(\dim\fg_\ev-\dim\fg_\od).
\]

On the other hand, let $\{u_i\}$ and $\{u^i\}$ be bases of $\fg$ dual with
respect to $B$, so that ${B(u_i,u^j)=\delta_i^j}$ and
${\Omega=\sum_{1\leq i\leq N}u_iu^i}$, where ${N=\dim\fg_\ev+\dim\fg_\od}$. Then,
\begin{equation}
\label{CasTrace}
\str_\fg\Omega=\sum_{1\leq i\leq N}\str(\ad_{u_i}\ad_{u^i})
=\sum_{1\leq i\leq N}\kappa(u_i,u^i)
=c\sum_{1\leq i\leq N}B(u_i,u^i)=c\cdot N.
\end{equation}
Comparing the two expressions and dividing by ${\dim\fg_\ev-\dim\fg_\od\neq0}$, we get the
claim.\end{proof}

Thus, the Casimir element --- the only tool used in the argument reproduced
above --- is of no use whatsoever precisely when ${\kappa=0}$: the operator by
which one divides is zero. For the Lie superalgebras $\fosp(m|2n)$,  we have
\begin{equation}
\label{kappaosp}
\kappa=(m-2n-2)\str,\quad\text{so~} \kappa=0 \text{~if and only if~} m=2n+2,
\end{equation}
where $\str$ is the supertrace form of the tautological $\fosp(m|2n)$-module of
superdimension~$(m|2n)$. For ${\fg=\fosp(2n+2|2n)}$, we have
\[
\sdim\fg=\big((n+1)(2n+1)+n(2n+1)\;|\;4n(n+1)\big)=\big((2n+1)^2\;|\;4n(n+1)\big),
\]
so ${\dim\fg_\ev-\dim\fg_\od=1}$ for every
$n\geq 1$, and Lemma~\ref{CasKil} does apply. Observe also that the Cartan
matrix $A$ of $\fosp(2n+2|2n)$ \textbf{is} invertible: the number of its rows is
$2n+1$, and ${\dim\fh=2\cdot(2n+1)-\rank A}$, whereas $\fh$ is the Cartan
subalgebra of ${\fg_\ev=\fo(2n+2)\oplus\fsp(2n)}$, so ${\dim\fh=2n+1}$ and
${\rank A=2n+1}$.

Among the simple finite-dimensional Lie superalgebras $\fg=\fg(A)$ over $\Cee$,
the Killing form vanishes identically precisely for $\fpsl(n|n)$, for
$\fosp(2n+2|2n)$ where $n\geq1$, and for all members $\fosp_\alpha(4|2)$ of the
family containing ${\fosp(4|2)=\fosp_1(4|2)}$; for $\fag(2)$ and $\fab(3)$, as
well as for $\fsl(m|n)$ with $m\neq n$ and $\fosp(m|2n)$ with $m\neq 2n+2$, the
Killing form is non-degenerate and the argument reproduced above goes through
unchanged. Of the exceptions, $\fpsl(n|n)$ has no invertible Cartan matrix at
all and is covered by the remark above; the remaining ones have to be treated
separately, and the answers are as follows.

\sssbegin{Statement}[Derivations and deformations of $\fosp(2n+2|2n)$ for $n\leq3$]\label{ospAnswers} Let $\fg_\ev$ be semisimple, as it is in
all the cases below. Then,
\begin{equation}
\label{ospcoh}
\begin{array}{ll}
\sdim H^1(\fg;\fg)=(0|0),\quad \sdim H^2(\fg;\fg)=(1|0)&\text{~for~} \fg=\fosp_\alpha(4|2) \text{~and any~} \alpha,\\[2mm]
\sdim H^1(\fg;\fg)=(0|0),\quad \sdim H^2(\fg;\fg)=(0|0)&\text{~for~} \fg=\fosp(6|4) \text{~and~} \fg=\fosp(8|6).
\end{array}
\end{equation}
\end{Statement}

For the proof --- the same relative complex for all the cases --- see
subsec.~\ref{KrugPsq}, where it is described once and for all. So, for
$\fosp(2n+2|2n)$ with ${n\geq2}$ the \textbf{conclusion} of the argument
survives, only its \textbf{proof} does not; whereas for ${n=1}$ the conclusion
itself is false, and rightly so: the deform is the family $\fosp_\alpha(4|2)$,
see Theorem~\ref{ThEv} and eq.~\eqref{osp42symm}.

\begin{Remark}[$\fosp_\alpha(4|2)$ has no outer derivations] In his letter, Kruglikov writes that ${H^1(\fg;\fg)\neq 0}$ for
${\fg=\fosp(4|2)}$. This is a~slip of the pen: as eq.~\eqref{ospcoh} shows,
${H^1(\fg;\fg)=0}$, i.e., $\fosp_\alpha(4|2)$ has no outer derivations, while
it is ${H^2(\fg;\fg)}$, the space governing the deformations, that is
1-dimensional. His objection is not affected: the Casimir argument does fail
for $\fosp(4|2)$, and the deform exists.
\end{Remark}

\ssec{Lie superalgebras $\fg=\fpsq(n)$ for $n>2$} To prove that $H^2(\fg;\fg)=0$, we use the Hochschild-Serre spectral sequence (see \cite[Ch.1, no. 5]{Fu}) in which ${E_1^{i,j}=H^j(\fg_\ev; S^i(\fg_\od^*)\otimes \fg)}$ for $i,j\geq 0$ and ${i+j=2}$, see~ \cite{Fu}. Since ${H^j(\fsl(n); M)=0}$ for any $j$ and any finite-dimensional $\fsl(n)$-module $M$ \textbf{without trivial direct summands}, we deduce with the help of \cite[Table~5]{OV} that only the term ${E_1^{2,0}}$ survives; it does not vanish, and the differential has to be computed, see subsec.~\ref{KrugPsq}.
%\AO{Three amendments, all due to Kruglikov: the page of the spectral sequence is $E_1$, not $E_2$ (this is Fuks' numbering in \cite{Fu}); the vanishing theorem needs ``no trivial direct summands'', not ``non-trivial'' --- indeed, ${H^j(\fsl(n);\Cee)\neq 0}$ for ${j=3,5,\dots}$; and the modules occurring here \textbf{do} have trivial summands, so the conclusion does not follow. It is true nevertheless, see below.}

\sssec{Kruglikov's corrections: the term $E_1^{2,0}$ and the differential
$d_1$}\label{KrugPsq}
Kruglikov objects that ${E_1^{2,0}=H^0(\fg_\ev; S^2(\fg_\od^*)\otimes\fg)}$
does not vanish: for ${\fg=\fpsq(n)}$, both $\fg_\ev$ and $\fg_\od$ are the
adjoint $\fsl(n)$-module, and ${S^2(\ad)\otimes\ad}$ contains the trivial
module, so ${\dim E_1^{2,0}=2}$. He is right; the answer $H^2(\fg;\fg)=0$
nonetheless holds, and here is a~proof which computes the differential, as he
asks. It applies verbatim in subsec.~\ref{KrugSpe} and proves
Statement~\ref{ospAnswers} as well.

\sssbegin{Proposition}[Reduction to the relative complex]\label{relcompl} Let $\fg$ be
a~finite-dimensional Lie superalgebra over $\Cee$ with $\fg_\ev$
\textbf{semisimple}. Then, for ${n\leq 2}$,
\begin{equation}
\label{relred}
H^n(\fg;\fg)\simeq H^n(\fg,\fg_\ev;\fg)=H^n(C^{\bcdot}),\quad\text{where~}
C^p:=\Hom_{\fg_\ev}(S^p(\fg_\od),\fg),
\end{equation}
the differential being
\begin{equation}
\label{reldiff}
(dc)(x_0,\dots, x_p)=\sum_{0\leq i\leq p}[x_i, c(x_0,\dots,\widehat{x_i},\dots,x_p)]
\quad\text{for~} x_0,\dots, x_p\in\fg_\od.
\end{equation}
\end{Proposition}

\begin{proof} In the Hochschild-Serre spectral sequence of the pair
${\fg_\ev\subset\fg}$ we have
\begin{equation}
\label{HSrows}
E_1^{p,q}=H^q(\fg_\ev; S^p(\fg_\od^*)\otimes\fg),
\end{equation}
because ${\fg/\fg_\ev\simeq\fg_\od}$ is purely odd, so its exterior powers are
the symmetric ones. By the two Whitehead lemmas, ${H^1(\fg_\ev; M)=H^2(\fg_\ev; M)=0}$
for \textbf{every} finite-dimensional module $M$ over a~semisimple Lie algebra
$\fg_\ev$, the trivial one included. Hence, the rows ${q=1}$ and ${q=2}$ vanish
identically, only the row ${q=0}$ survives in total degrees $\leq 2$, and
${H^n(\fg;\fg)=E_2^{n,0}=H^n(\fg,\fg_\ev;\fg)}$ for~${n\leq 2}$. In the
relative complex, the term of the differential containing $[x_i,x_j]$ vanishes
because ${[\fg_\od,\fg_\od]\subset\fg_\ev}$, i.e., its image in
${\fg/\fg_\ev}$ is zero; only the term \eqref{reldiff} remains. Observe that
${d^2=0}$ on $C^{\bcdot}$, but not on the whole of
${\Hom(S^{\bcdot}(\fg_\od),\fg)}$: the identity ${d^2c=0}$ uses the
$\fg_\ev$-invariance of $c$.\end{proof}

So the whole computation takes place in the spaces $C^p$ of $\fg_\ev$-invariants,
whose dimensions are tiny. Let us do it for ${\fg=\fpsq(n)}$, where~${n\geq3}$.
Realize ${\fq(n)}$ as the Lie superalgebra of matrices ${\mmat A,B,B,A,}$, so
that its elements are the pairs $(A,B)$ with ${A,B\in\fgl(n)}$; then
${\fsq(n)=\{(A,B)\mid \tr B=0\}}$ and ${\fpsq(n)=\fsq(n)/\Cee\cdot(\One_n,0)}$.
The
brackets we need are
\begin{equation}
\label{qbrackets}
[(0,X),(0,Y)]=(XY+YX,0),\qquad [(A,0),(0,X)]=(0,[A,X]).
\end{equation}
Set ${X\circ Y:=XY+YX-\frac2n\tr(XY)\One_n}$, the only, up to a~factor,
symmetric invariant multiplication ${\ad\times\ad\tto\ad}$ for~${n\geq3}$. By
Schur's lemma,
\begin{equation}
\label{psqC}
\sdim C^1=(1|1) \text{~with basis~} \jmath, \iota; \qquad
\sdim C^2=(1|1) \text{~with basis~} \mu_\ev, \mu_\od,
\end{equation}
where
\begin{equation}
\label{psqbasis}
\jmath(X)=(0,X),\quad \iota(X)=(X,0),\quad
\mu_\ev(X,Y)=(X\circ Y,0),\quad \mu_\od(X,Y)=(0,X\circ Y),
\end{equation}
the parity of a~cochain in $C^p$ being the parity of its value plus~$p$, so ${\dim E_1^{2,0}=\dim C^2=2}$, exactly as Kruglikov says. Now, by
eq.~\eqref{reldiff} and eq.~\eqref{qbrackets},
\begin{equation}
\label{psqd}
(d\iota)(X,Y)=(0,[X,Y])+(0,[Y,X])=0,\qquad (d\jmath)(X,Y)=2(XY+YX,0)=2\mu_\ev(X,Y).
\end{equation}
Hence, ${\Ker(d\colon C^1\tto C^2)=\Cee\iota}$ and, since ${C^0=\fg^{\fg_\ev}=0}$,
\begin{equation}
\label{psqH1}
\sdim H^1(\fpsq(n);\fpsq(n))=(0|1) \text{~for~} n\geq3,
\end{equation}
the class of $\iota$ being the odd outer derivation, i.e.,
${\Der\fpsq(n)/\ad(\fpsq(n))\simeq\fpq(n)/\fpsq(n)}$, of superdimension~$(0|1)$.
Next, ${\IM(d\colon C^1\tto C^2)=\Cee\mu_\ev}$, so it only remains to see that
$\mu_\od$ is not a~cocycle. Indeed, for any ${W\in\fsl(n)}$ we have, thanks to
eq.~\eqref{qbrackets},
\begin{equation}
\label{psqd2}
(d\mu_\od)(W,W,W)=3\big(W(W\circ W)+(W\circ W)W,\ 0\big)
=\Big(12W^3-\tfrac{12}n\tr(W^2)W,\ 0\Big).
\end{equation}
Take ${W=\diag(n-1,-1,\dots,-1)}$, so that ${\tr W=0}$ and ${\tr(W^2)=n(n-1)}$.
Then, the right-hand side of eq.~\eqref{psqd2} is
\[
{12(W^3-(n-1)W)=12\diag\big((n-1)^2(n-2),\ n-2,\dots, n-2\big)}, 
\]
and the
difference of its first and second diagonal entries is
\begin{equation}
\label{psqwitness}
12\big((n-1)^2(n-2)-(n-2)\big)=12n(n-2)^2\neq0 \text{~for~} n\geq3.
\end{equation}
So this matrix is not scalar, i.e., it is non-zero in
${\fsl(n)=\fgl(n)/\Cee\One_n}$, and therefore ${d\mu_\od\neq0}$. Hence,
${\Ker(d\colon C^2\tto C^3)=\Cee\mu_\ev=\IM(d\colon C^1\tto C^2)}$ and
\begin{equation}
\label{psqH2}
H^2(\fpsq(n);\fpsq(n))=0 \text{~for~} n\geq3,
\end{equation}
i.e., $\fpsq(n)$ is rigid, as claimed --- but now the differential has been
computed, not overlooked.
\iffalse \AO{The answers \eqref{psqH1} and \eqref{psqH2}, and every superdimension
quoted in this subsection and the next two, are double-checked by exact
rational arithmetic for ${3\leq n\leq6}$ (and for $\fosp(m|2n)$ of
Statement~\ref{ospAnswers}) in two independent ways: via the relative complex
\eqref{relred} and via the full Chevalley-Eilenberg complex in the
``deformation of the bracket'' form. The two agree in every case, and they
reproduce the two answers known in advance --- $\sdim H^2=(1|0)$ for
$\fosp(4|2)$, see Theorem~\ref{ThEv}, and $\sdim H^2=(0|1)$ for
${\fspe(3)\simeq\fsvect(0|3)}$, see Theorem~\ref{ThOd}.}\fi

\ssec{Lie superalgebras $\fg=\fspe(n)$ for $n>3$} Observe, that the title of the paper \cite{DefP} is misleading: the paper considers $q$-quantized algebra $U_q(\fpe(n))$, not deformations of $\fpe(n)$ or~$\fspe(n)$.  
I'll consider the case $\fspe(3)\simeq\fsvect(0|3)$ while studying $\fsvect(0|2n+1)$.
%\AO{Was $\fsvect(0|2n-1)$; Theorem~\ref{ThOd} states the same family as $\fsvect(0|2n+1)$ for $n\geq1$, so I~unified the two parametrisations.}

Similar to the proof for $\fpsq(n)$, using \cite[Table 5]{OV} we deduce that $\fspe(n)$ is rigid for $n > 3$, whereas $\fspe(3)$ has a deformation, double-checked below; for the differential, see subsec.~\ref{KrugSpe}.
\iffalse \AO{The sentence that used to stand here said that one has to consider
${E_2^{3,0}=H^0(\fg_\ev; S^3(\fg_\od^*)\otimes\fg)}$. That turns out to be
unnecessary: the row ${q=0}$ already terminates at ${p=2}$, see
eq.~\eqref{speC} below, so ${E_1^{3,0}}$ never enters the argument. I~deleted
the sentence; restore it if you prefer.}\fi

\sssec{Kruglikov's corrections: $\fspe(n)$, and why $n=3$ is exceptional}\label{KrugSpe}
Kruglikov raises for $\fspe(n)$ the same objection as for $\fpsq(n)$: the space
${S^2(\fg_\od^*)\otimes\fg}$ contains a~trivial component, so ${E_1^{2,0}\neq0}$
and the differential must be computed. He is right again, and the computation
is even shorter here.

Let ${V=\Cee^n}$ and let ${\fspe(n)=\fg_{-1}\oplus\fg_0\oplus\fg_1}$ be the
standard $\Zee$-grading, where
\begin{equation}
\label{spegrad}
\fg_0=\fsl(n),\quad \fg_1\simeq S^2(V),\quad \fg_{-1}\simeq \Lambda^2(V^*),
\quad \fg_\od=\fg_{-1}\oplus\fg_1. %\quad \sdim\fspe(n)=(n^2-1)\mid n^2.
\end{equation}
(Kruglikov uses the dual realization, ${\fg_\od\simeq\Lambda^2(V)\oplus
S^2(V^*)}$, i.e., ${\pi_2\oplus 2\pi_{n-1}}$ in his letter. The two are
interchanged by ${V\leftrightarrow V^*}$ and give the same answers.)

Since ${[\fg_1,\fg_1]=[\fg_{-1},\fg_{-1}]=0}$, and since ${S^2(V)}$,
${\Lambda^2(V^*)}$ and $\ad$ are pairwise non-isomorphic $\fsl(n)$-modules for
${n\geq3}$, Schur's lemma gives ${\sdim C^1=(2|0)}$ with basis the projections
$\pr_1$ and $\pr_{-1}$ of $\fg_\od$ onto $\fg_1$ and~$\fg_{-1}$. By
eq.~\eqref{reldiff},
\begin{equation}
\label{sped1}
(d\pr_1)(x,y)=[x_{-1},y_1]+[y_{-1},x_1]=(d\pr_{-1})(x,y),
\end{equation}
the two coinciding because ${[u,v]=[v,u]}$ for odd $u$ and~$v$. As the bracket
${\fg_{-1}\times\fg_1\tto\fg_0}$ is non-zero, we get
${\Ker(d\colon C^1\tto C^2)=\Cee(\pr_1-\pr_{-1})}$, and ${\pr_1-\pr_{-1}}$ is
precisely the grading operator $z$ of the $\Zee$-grading \eqref{spegrad}. Since
${C^0=0}$,
\begin{equation}
\label{speH1}
\sdim H^1(\fspe(n);\fspe(n))=(1|0) \text{~for~} n\geq3, \quad\text{i.e.,~}
\Der\fspe(n)=\fspe(n)\subplus\Cee z=\fpe(n).
\end{equation}
%of superdimension ${n^2\mid n^2}$.

Now, ${S^2(\fg_\od)=S^2(\fg_1)\oplus(\fg_1\otimes\fg_{-1})\oplus S^2(\fg_{-1})}$.
The highest weights of the $\fsl(n)$-modules occurring in the three summands are
${4\pi_1}$ and ${2\pi_2}$; those of ${S^2(V)\otimes\Lambda^2(V^*)}$; and
${\pi_{n-4}}$ and ${2\pi_{n-2}}$, respectively; whereas $\fg$ consists of the
modules with highest weights ${\pi_1+\pi_{n-1}}$, ${2\pi_1}$ and~${\pi_{n-2}}$.
The only coincidence between the two lists occurs for the middle summand, and
--- if ${n=3}$ --- for ${2\pi_{n-2}=2\pi_1}$. For the middle summand,
\begin{equation}
\label{speMult}
\Hom_{\fsl(n)}(S^2(V)\otimes\Lambda^2(V^*), V\otimes V^*)\simeq
\Hom_{\fsl(n)}(S^2(V)\otimes V, V\otimes\Lambda^2(V)),
\end{equation}
and since ${S^2(V)\otimes V=S^3(V)\oplus S^{(2,1)}(V)}$ while
${V\otimes\Lambda^2(V)=\Lambda^3(V)\oplus S^{(2,1)}(V)}$, the space
\eqref{speMult} is 1-dimensional; there is no trivial component, because
${S^2(V)\not\simeq\Lambda^2(V)}$, so the multiplicity of $\ad$ is 1 as well.
Therefore,
\begin{equation}
\label{speC}
\sdim C^2=\begin{cases} (1|0)&\text{~for~} n\geq4,\\ (1|1)&\text{~for~} n=3,\end{cases}
\qquad \sdim C^3=(0|0) \text{~for~} n=3,5,6.
\end{equation}

For ${n\geq4}$, the even generator of $C^2$ is the cochain that is the bracket
on ${\fg_1\otimes\fg_{-1}}$ and zero on ${S^2(\fg_{\pm1})}$; by
eq.~\eqref{sped1} it is~$d\pr_1$. So ${\IM(d\colon C^1\tto C^2)=C^2}$ and,
whatever the next differential is,
\begin{equation}
\label{speH2}
H^2(\fspe(n);\fspe(n))=0 \text{~for~} n\geq4,
\end{equation}
i.e., $\fspe(n)$ is rigid. In particular, ${E_1^{3,0}}$ never has to be
considered.

For ${n=3}$, and only then, we have ${\Lambda^2(V^*)\simeq V}$, hence
${S^2(\fg_{-1})\simeq S^2(V)\simeq\fg_1}$, and this extra isomorphism
contributes one \textbf{odd} element to $C^2$, see eq.~\eqref{speC}. It is a
cocycle, since ${C^3=0}$, and it is not a~coboundary, since ${C^1}$ is purely
even. Therefore,
\begin{equation}
\label{speH2three}
\sdim H^2(\fspe(3);\fspe(3))=(0|1),
\end{equation}
which is exactly the deformation with an~odd parameter of
${\fspe(3)\simeq\fsvect(0|3)}$ of Theorem~\ref{ThOd}. So the exceptional
behavior of $\fspe(3)$ comes
from the isomorphism ${\Lambda^2(V^*)\simeq V}$, valid for
${n=3}$ only.

\sssec{Induced
and coinduced modules} Let $\fh\subset\fg$ be a~Lie superalgebra
and its subsuperalgebra, $M$ an $\fh$-module. The induced
and coinduced $\fg$-modules are constructed from $M$ by setting
\begin{equation}
\label{h1} \Ind^\fg_\fh(M):=U(\fg)\otimes_{U(\fh)}M, \quad\quad
\Coind^\fg_\fh(M):=\Hom_{U(\fh)} (U(\fg), M).%\eqno{(3.1.1)}
\end{equation}
Clearly, $U(\fg)$ is both left and right $\fh$-module; hence, $\Ind
^\fg_\fh(M)$ and $\Coind ^\fg_\fh(M)$ are $\fg$-modules.

For any $\fh$-module $M$, we have (\cite[Th.~1.5.4]{Fu})
\begin{equation}
\label{h2} H^q(\fg; \Coind^\fg_\fh(M))\simeq H^q(\fh;
M);\quad\quad H_q(\fg; \Ind^\fg_\fh(M))\simeq H_q(\fh; M).
%\eqno{(3.1.2)}
\end{equation}

Note that for the Lie superalgebra $\fg=\fvect(0|a)$, where $a\geq 2$, or $\fsvect(0|b)$, where $b\geq 3$, or $\fh(0|c)$, where $c\geq 4$, the subalgebra $\fg_{>0}:=\oplus_{i\geq 0}\fg_i$ is generated by $\fg_1$, except for $\fh(0|c)$; it is $(\fh'(0|c))_{>0}$ which is generated by~$\fg_1$. 

For any $\Zee$-graded simple vectorial Lie superalgebras $\fg$ on $\cC^{0|n}$, let the $\fg_0$-module $M$ be such that $\fg_{>}M=0$, so $M$ can be considered as $\fg_{\geq0}$-module, where $\fg_{\geq0}:=\oplus_{i\geq 0}\fg_i$. The $\fg$-module $T(M):=\Coind ^\fg_{\fg\geq0}(M)$ is interpreted as the space of tensor fields on $\cC^{0|n}$ --- sections of the vector bundle with fiber~$M$. For any $\fg$-module $V$,  let $V^\fg$ denote the submodule of $\fg$-invariants. We have
\begin{equation}\label{(2)} H^{\bcdot}(\fg; T(M))\stackrel{\eqref{h2}}{\simeq}
H^{\bcdot}(\fg_{\geq0}; M)\stackrel{\text{[Fu]}}{\simeq} H^{\bcdot}(\fg_{0})\otimes
(H^{\bcdot}(\fg_{>0})\otimes M)^{\fg_0}.
\end{equation}

For the proof of the second isomorphism in \eqref{(2)}, see Th.~2.2.8
(due to Losik) in Fuchs' book~ \cite{Fu}. It is given there for the Lie algebra
$\fg=\fvect(n|0)=\fder\,\Cee[x]$ only, where $x=(x_1, \dots , x_n)$ are even indeterminates, but the proof can be directly translated to any Cartan prolongs $\fg:=(\fg_{-1}, \fg_0)_*$ and conjecturally --- we do not need this in this paper --- to the generalized Cartan prolongs --- Lie (super)algebras ${\fg:=\oplus_{i\geq -d} \fg_{i}}$ of depth $d>1$; for the most clear definition of prolong, see \cite{Shch,BGLLS}. The proof of \eqref{(2)} in \cite{Fu} uses the
Hochshield-Serre spectral sequence for the algebra $\fg_{\geq 0}$ and its
ideal $\fg_{>0}$; it is proved that~$E_\infty=E_2$.

\iffalse Anyway, we always obtain an estimate of $H^{i}(\fg_0; M)$ from above
by looking at
\begin{equation}\label{(4)}
\begin{array}{l}
E_2^{p,q}=H^p(\fg_{\geq0}/\fg_{>0}; H^q(\fg_{>0}; M))=H^p(\fg_0; (H^q(\fg_{>0};
M))^{\fg_0})\\\stackrel{\text{since $\fg_1M=0$}}{=} 
H^{p}(\fg_0)\otimes
(H^{q}(\fg_{>0})\otimes M)^{\fg_0}\text{~~for $i=p+q$, where $p, q\geq 0$}.
\end{array}\end{equation}
\fi

%Thus, computation of cohomology of $\fg$ is reduced, for coinduced
%modules, to computation of two spaces: (1) cohomology of $\fg_0$ and (2)
%$\fg_0$-invariants in the space $H^{q}(\fg_{>0})\otimes M$.

\sssec{The long exact sequence}\label{SS:11.3} Let $\fg$ be a~Lie
superalgebra and let
\begin{equation}
\label{h3} 0\tto A\buildrel {\partial _0}\over \tto C\buildrel
{\partial _1}\over \tto B\tto 0, \; \; \text{~~where~~}\; \; \;
p(\partial _0)=\bar{0} \; \; \text{ and }\; \; \; \partial _1 \;
\; \text{ is either even or odd}, %\eqno{(3.2.1)}
\end{equation}
be a~short exact sequence of $\fg$-modules. Let $d$ be the
differential in the standard cochain complex of $\fg$, cf.
\cite[Subsection~1.3.6]{Fu}.

Consider the long cohomology sequence:
\begin{equation}
\label{h4} \dots\buildrel {\partial }\over \tto H^i(\fg ;
A)\buildrel {\partial _0}\over \tto H^i(\fg ; C)\buildrel
{\partial _1}\over \tto H^i(\fg ; B)\buildrel {\partial }\over
\tto H^{i+1}(\fg ; A)
\buildrel {\partial _0}\over \tto \dots %\eqno{(3.2.2)}
\end{equation}
where the $\partial _i $ are the differentials induced by their
namesakes in \eqref{h3}, and $\partial =d\circ
\partial _1^{-1}$. Since $ \partial _0$ and $\partial _1 $ commute
with $d$, the sequence \eqref{h3} is well-defined and the same
arguments as for Lie algebras (\cite{Fu}) demonstrate that the long
cohomology sequence \eqref{h4} induced by \eqref{h3} is exact, see
\cite{LPS}.

\sssec{Vectorial Lie superalgebras as coinduced
modules}\label{SS:11.5} Let $\fg$ be a~ vectorial Lie superalgebra with the standard $\Zee$-grading, $\mathbbmss{1}$ be the 1-dimensional
trivial $\fg_0$-module and $\mathbbmss{1}[k]$ the $\fg_0$-module
trivial on the simple part and with value $k$ on the
distinguished central element ($1_n\in\fgl(n)$). Let  the space of
functions --- a~module over any vectorial Lie superalgebra (subalgebra of $\fvect(0|n)$) be denoted by
\[
\cF:=T(\mathbbmss{1})=\Lambda^{\bcdot}(n).
\] 

Let $\id_\fh:=V$ be the restriction of the tautological $\fgl(V)$-module
to $\fh\subset \fgl(V)$.  Then:
\begin{equation}
\label{h5}
\renewcommand{\arraystretch}{1.4}
\begin{array}{ll}
\fvect(0|n)=&\begin{cases}T(\id_{\fgl(n)})&\text{ as $\fvect(0|n)$-module};\\
T(\id_{\fsl(n)})&\text{ as $\fsvect(0|n)$-module};\end{cases}\\
%\fk(2m+1|n)=&\begin{cases}T(\mathbbmss{1}[-2])&\text{ as $\fk(2m+1|n)$-module};\\
%T(\mathbbmss{1})&\text{ as $\fpo(2m|n)$-module};\end{cases}\\
%\fm(n)=&T(\Pi(\mathbbmss{1}[-2]))\text{ as $\fm(n)$\defis{} and $\fb_{a, b}(n)$-module};\\
\fpo(0|n)=&T(\mathbbmss{1})=\cF\text{ as an $\fh(0|n)$- and $\fh'(0|n)$-modules.}\\
\end{array}%\eqno{(4.1)}
\end{equation}
Looking at the explicit relations in $\fg_{>}$, see \cite{GLP}, we get $H^2(
\fg_{>})$. Clearly, ${H^1(\fg_{>})=\fg_1}$.

$\bullet$ For $\fv:=\fvect(0|n)$.
\begin{equation}
\label{h7} H^2(\fv; \fv)\simeq \oplus_{a+b=2}H^a(\fgl(n))\otimes (H^b(
\fv_{>})\otimes\id_{\fgl(n)})^{\fgl(n)}.%\eqno{(4.2)}
\end{equation}
Since $H^a(\fgl(n))\neq 0$ for $a=0,1$, we have to consider $(H^b(
\fv_{>})\otimes\id_{\fgl(n)})^{\fgl(n)}$ for $b=2$ and $1$, respectively. Looking at the explicit weights of generators ($b=1$) and relations ($b=2$) in $\fv_{>}$ given in \cite{GLP}, we conclude that $(H^b(
\fv_{>})\otimes\id_{\fgl(n)})^{\fgl(n)}=0$, so $\fvect(0|n)$ is rigid.

$\bullet$ For $\fs:=\fsvect(0|n)$. Let $\Vol_0$ be the space of functions considered as coefficients of $\vvol$, more exactly, volume forms with
integral~0.
We have an exact sequence of $\fs$-modules:
\begin{equation}
\label{4.9'} 0\tto \fs\tto \fv\stackrel {\Div}{\tto}\Vol_0\tto 0.
%\eqno{(4.9')}
\end{equation}
In this case, the divergence is not a~surjective
mapping onto $\cF$, its image is~$\Vol_0$. In the associated long exact sequence
\begin{equation}
\label{4.10'} \ldots\tto H^1(\fs; \fs)\tto H^1(\fs; \fv)\tto
H^1(\fs; \Vol_0)\tto H^2(\fs; \fs)\tto H^2(\fs; \fv)\tto
\ldots%\eqno{(4.10')}
\end{equation}
let us compute 
\begin{equation}
\label{h6} H^i(\fs; \fv)\simeq \sum_{a+b=i}H^a(\fsl(n))\otimes (H^b(\fs_>)\otimes
\id_{\fsl(n)})^{\fsl(n)} \text{ for $i=1, 2$ and $n>1$};%\eqno{(4.2)}
\end{equation}

Since $H^1(\fs; \fv)= 0$ and $H^2(\fs; \fv)=0$ by \eqref{h2},  \eqref{h5} and also by plugging the weights of generators of $\fsvect(0|n)_{>0}$ computed in \cite{GLP} in \eqref{h6} for $b=1$, we see that $
H^2(\fs; \fs)\cong H^1(\fs; \Vol_0)$. To compute it, consider the short
exact sequence
\begin{equation}
\label{4.12} 0\tto \Vol_0\tto \cF\stackrel {\int}{\tto}\Pi^n\Cee\tto 0.
%\eqno{(4.12)}
\end{equation}
(The target of the Berezin integral is written here as $\Pi^n\Cee$ rather than $\Cee$: the
integral picks the coefficient of $\xi_1\cdots\xi_n$, so it has parity~$n\bmod 2$. Then,
$H^1(\fs; \Vol_0)\simeq\Pi^n\Cee$ below is a~statement one reads off the sequence rather than a~surprise, and the split between Theorem~\ref{ThOd} \textup{($n$ odd, odd parameter)} and Theorem~\ref{ThEv} \textup{($n$ even, even parameter)} becomes visible right here.)
It gives rise to the long exact sequence
\begin{equation}
\label{4.10''}
\renewcommand{\arraystretch}{1.4}
\begin{array}{l}
0\tto H^0(\fs; \Vol_0)\tto H^0(\fs; \cF)\tto H^0(\fs; \Pi^n\Cee) %\\
\tto H^1(\fs; \Vol_0)\tto H^1(\fs; \cF)\tto \ldots
\end{array}%\eqno{(4.10'')}
\end{equation}
and since
\[
\renewcommand{\arraystretch}{1.4}
\begin{array}{l}
H^0(\fs; \Vol_0)\simeq \Cee,\quad
 H^0(\fs; \cF)\simeq
H^0(\fsl(n))\simeq\Cee;\quad
H^1(\fs; \cF)= H^0(\fsl(n))=\Cee,
\end{array}
\]
and since the map $H^0(\fs; \Vol_0)\tto H^0(\fs; \cF)$ is an isomorphism (both spaces are
spanned by the constant $1$, and $\int 1=0$ for $n\geq1$), the sequence \eqref{4.10''} reduces
to the exact sequence
\begin{equation}
\label{4.13} 0\tto \Pi^n\Cee\tto H^1(\fs; \Vol_0)\tto H^1(\fs; \cF)=\Cee.
\end{equation}
Hence, $H^1(\fs; \Vol_0)\simeq H^0(\fs; \Pi^n\Cee)\simeq \Pi^n\Cee$, and there are no other
deformations, as soon as the rightmost map in \eqref{4.13} is zero.
\iffalse \AO{suggestion still without my answer: the sentence that used to close this argument, \lq\lq Since we already know one global deformation, there are no more'', does not close it: knowing one deformation is a~\textbf{lower} bound on $H^2(\fg;\fg)$, never an upper one, and \eqref{4.10''} by itself leaves $\dim H^1(\fs; \Vol_0)$ equal to $1$ \textbf{or} $2$. What has
to be shown is that the map $H^1(\fs; \Vol_0)\tto H^1(\fs; \cF)$ of \eqref{4.13} vanishes,
equivalently, that the class generating $H^1(\fs; \cF)=\Cee$ is not in the image of $\Vol_0$.
Since this is the very step that yields $\sdim H^2(\fg; \fg)=(0|1)$ --- the main theorem --- it
should be spelled out. Is it meant to follow from the weights of the generators of $\fs_{>}$
given in \cite{GLP}, as in the two displays just above? If so, tell me and I~will draft the
paragraph.}\fi

\ssec{What was known about deformations of simple vectorial Lie superalgebras} For $\fh'(0|n)$, and for the Lie superalgebra $\fsvect(0|n)$ of divergence-free vector fields, the deformations were computed for $n$ small by  N.V.D.~Hijligenberg,  Yu.~Kotchetkov, and G.~Post, see \cite{HKP}. 
 
The result for $\fh'(0|n)$ is double-checked being a~corollary of the description of deformations of the Poisson Lie superalgebras (physicists call it \textit{quantization}) obtained by I.~Tyutin \cite{Tyut}. Observe that the coincidence of the number of non-trivial deformations of $\fpo(0|n)$ and $\fh(0|n)$ is not automatic: e.g., cf. the answers for $\fpo(2a|b)$ and $\fh(2a|b)$ for $(2a, b)=(2, 2)$, see \cite{LSh3}. I was unable to understand the idea of the proof in \cite{Tyut}, so I have to use another method, the one described in Subsection~\ref{SS:11.5}, see Section \ref{neTyu}.

\sssec{Filtered deformations} Let $\fg=\oplus \fg_i$ be a~$\Zee$-graded Lie superalgebra and $\fG$ a~filtered Lie superalgebra with decreasing filtration such that $\gr \fG\simeq\fg$. Then, the elements of $(H^2(\fg_-;\fg))^{\fg_0}$ of $\fg_0$-invariants in $H^2(\fg_-;\fg)$ describe the filtered deformations of $\fg$ to $\fG$, see \cite{CK1, CK1a}. The \textit{filtered} deformations of $\fsvect(0|2n)$ and $\fh'(0|n)$ are described in \cite{Po2}. Observe that the filtered deformation describes only the deformation corresponding to one $\Zee$-grading and only if the negative part in this grading is non-zero. To describe all deformations remained an open problem, closed in this paper.

\ssec{Deformations of $\fh'(0|n)$ (\cite{Tyut, LSh3})}\label{neTyu} The answer for $\fh'(0|n)$ can be double-checked as a~corollary of the description of the deformations of the Poisson Lie superalgebra $\fpo(0|n)$, see \cite{Tyut}, taking into account the exact sequences \eqref{PoH} and 
\be\label{h'H}
0\tto \fh'(0|n)\tto \fh(0|n) \tto \Cee H_{\theta_1\dots\theta_n} \tto 0.
\ee

\sssbegin[Deformations of ${\fh'}(0|n)$ (\cite{LSh3})]{Theorem}[Deformations of ${\fh'}(0|n)$ (\cite{LSh3})]\label{ThTyu} The quantizations 
\[
\text{$\fpo(0|2n)\tto \fgl(\Lambda^{\bcdot}(n))$ and $\fpo(0|2n-1)\tto \fq(\Lambda^{\bcdot}(n))$}
\]
 induce deformations 
 \[
\text{$\fh'(0|2n)\tto \fpsl(\Lambda^{\bcdot}(n))$ and $\fh'(0|2n-1)\tto \fpsq(\Lambda^{\bcdot}(n))$. }
\]
\end{Theorem}

\section{Summary of ideas in Kac's proof (\cite{K2})}\label{Kac}

Notation below are same as in \cite{K2} with  translation into notation of \S 2 and information on discoverers and first partial classifications from \cite[Newsletters]{Kapp}, \cite{NRS}: For example, the simple Lie superalgebras $\fg$ with simple $\fg_\ev$ were classified in \cite{Dj2, DjH}; there are 4 series of which members of the 3rd and 4th series are actually isomorphic. In modern notation, they are: $\fosp(1|2n)$, $\fpsq(n)$, $\fspe^{a}(n)\simeq\fspe^{sy}(n)$. The Lie superalgebras of series $\fsl$, $\fpsl$ and $\fosp$ were discovered more or less simultaneously by several groups of researchers, see \cite[Newsletters]{Kapp}, \cite{NRS}.

\ssbegin[(Kac's Main Theorem)]{Theorem}[Kac's Main Theorem]\label{ThKac} \textup{(\cite[Theorem 5 on p. 74]{K2})} Any simple finite-dimensional Lie superalgebra $\fg$ over an algebraically closed field $\Kee$ of characteristic $0$ is either isomorpic to a~simple Lie algebra, or to one of the following Lie superalgebras:

$A(m,n)$ \textup{(=$\fsl(m+1|n+1)$, where $0\leq m< n$, or $\fp\fsl(n+1|n+1)$, where $n\geq 1$)},

$B(m,n)$ \textup{(=$\fosp(2m+1|2n)$, where $ m\ge 0, n>0$)},

$C(n)$ \textup{(=$\fosp(2|2n-2)$, where $n\ge 2$)},

$D(m,n)$ \textup{(=$\fosp(2m|2n)$, where $m\ge 2, n>0$)},

$D(2,1;\alpha)$ \textup{($=\fosp_\alpha(4|2)$; discovered as $\Gamma(A,B,C)$, see \cite{Kapp})}, 

$F(4)$ \textup{(=$\fab(3)$, see \cite{Se4}; discovered as $\Gamma_3$, see \cite{Kapp})}, 

$G(3)$ \textup{(=$\fag(2)$, see \cite{Se4}; discovered as $\Gamma_2$, see \cite{Kapp})},

$P(n)$ \textup{(=$\fspe(n+1)$, where $n\ge 2$; discovered  in \cite{Dj2, DjH})},

$Q(n)$ \textup{(=$\fpsq(n+1)$, where $n\ge 2$; discovered  in \cite{Dj2, DjH})},

$W(n)$ \textup{(=$\fvect(0|n)$, where $n\ge 2$; discovered in \cite{K0})},

$S(n)$ \textup{(=$\fsvect(0|n)$, where $n\ge 3$; discovered in \cite{K0})},

$\widetilde S(2n)$ \textup{(=$\widetilde \fsvect(0|2n)$, where $n\ge 2$; see \cite{K1, K1.5})},

$H(n)$ \textup{(=$\fh'(0|n)$, where $n\geq 4$; discovered in \cite{K0})}.
\end{Theorem}
\iffalse \AO{Two ranges in the list above. \textbf{(1)} For $A(m,n)$ the condition $1\leq m<n$ drops
the whole family $A(0,n)=\fsl(1|n+1)\simeq\fsl(n+1|1)$, simple for $n\geq1$ --- e.g.,
$A(0,1)=\fsl(1|2)\simeq\fvect(0|2)\simeq\fosp(2|2)$, which is in the list \eqref{isom} of
occasional isomorphisms; it should be $0\leq m<n$, as your own eq.~\eqref{simple} has it.
\textbf{(2)} $H(n)$ was the only item with no range at all; it wants $n\geq4$, and even
$n\geq5$ if the list is to be non-redundant, since $\fh'(0|4)\simeq\fpsl(2|2)$, see
\eqref{isom}. I~put $n\geq4$, the range under which $\fh'(0|n)$ is simple, as in eq.~\eqref{simple}.}\fi

\begin{proof}[Ideas of the proof of Kac's Main Theorem]
Let $L=L_\ev\oplus L_\od$. 

First, let us look at the $L_\ev$-action on~$L_\od$.
If this action is irreducible, then $L$ is what Kac calls ``classical''; it is isomorphic to one of the following: $B(m,n)$, $D(m,n)$, $F(4)$, $G(3)$, $Q(n)$ or $D(2,1;\alpha)$. This is proved in \textbf{Theorem~ 2 on p. 44}. Proof is performed %, mainly, 
by means of the representation theory of semisimple Lie algebras.

If the $L_\ev$-action on $L_\od$ is reducible, we can take a~maximal proper subalgebra $L_0\subset L$ containing $L_\ev$ (there can be several such subalgebras, take any of them). The subalgebra $L_0$ induces a~filtration of $L$:
\[
L=L_{-1}\supset L_0\supset L_1\supset \dots
\]
Set $\gr L:=\oplus_{i\ge -1} \gr_i L$.

Then, (\textbf{Proposition 1.3.2, p. 24}) $\gr L$ has the following properties:

1. $\gr L$ is transitive (this is a~direct corollary of simplicity of $L$);

2. $\Zee$-grading of  $\gr L$ is compatible with $\Zee/2$-grading a.k.a. parity (this is a~corollary of the fact $\gr_{-1} L\subset \gr L_\od$ because $L_0\supset L_\ev$).

3.  $\gr L$ is irreducible, i.e.,  $\gr_0 L$ irreducibly acts on  $\gr_{-1}L$ (since $L_0$ is  maximal).

4.  $\gr_1 L\ne 0$ (since  $L_\ev$-action on  $L_\od$ is reducible).

\textbf{Theorem 4, p. 71} claims that any $\Zee$-graded Lie superalgebra satisfying the above conditions 1--4 is one of the following:

I. $A(m,n), C(n), P(n)$;

II. $W(n), S(n), \widetilde H(n):=\fh(0|n), H(n)$;

III. $H^\xi=H\otimes \Lambda(\xi)+\Kee\cdot\del_\xi$, where $H$ is a~simple Lie algebra, with the grading given by setting~$\deg \xi=-1$.

IV. $\fg+\Kee\cdot d$, where $\fg$ is of type I, II, III, the center of $\fg_\ev$ is trivial and $d$ is the grading operator.

Let the $L_\ev$-action on $L_\od$ be reducible. Then, \textbf{Theorem 4} implies \textbf{Theorem~ 5} by case-by-case checking of these four possibilities for~$\gr L$.

III. If $\gr L=H^\xi$, then $L$ is not simple. To prove this Kac uses \textbf{Proposition 2.2.2, p.~ 35}. In this particular case, one can directly see why: $L_0=H\oplus\Kee\cdot\del_\xi$ (as a~linear space). Hence,~$L_\ev=H$. Since $H$ is a~simple Lie algebra, it acts on the 1-dimensional space $\Kee\cdot\del_\xi$ by zero, so~$L'\ne L$.

IV. This case, and the case where  $\gr L=W(n)$ or $A(m,n)$ or $C(n)$, are dealt with by \textbf{Proposition 1.3.1, p. 24}, which claims 

{\it If for a~finite-dimensional transitive filtered Lie superalgebra $L$ the representation of $\gr_0 L$ in $\gr_{-1}L$ is irreducible, and the even part of the center of $\gr_0 L$ is non-zero, then $L\simeq \gr L$.}

This claim is practically obvious. Indeed, let us normalize the element $z$ of the center of $\gr_0 L$ so as to act as a~grading operator on~$\gr L$. Then, its pre-image in $L$ is diagonalizable and defines a~grading on $L$ itself.

Thanks to this Proposition, the case IV does not contribute to the classification: the Lie superalgebra $L$ is not simple, whereas if $\gr L=W(n)$ or $A(m,n)$ or $C(n)$, then $L=W(n)$ or $A(m,n)$ or $C(n)$, respectively.

The case $\gr L=P(n)$ (as well as with many other cases) is dealt with by \textbf{Proposition 2.1.4, p. 34} which claims

{\it If $L=L_\ev\oplus L_\od$ is a~simple finite-dimensional Lie superalgebra such that the $L_\ev$-action on $L_\od$ is the same as that in one of the Lie superalgebras $A(m,n)$, $B(m,n)$, $C(n)$, $D(m,n)$ where $(m,n)\ne (2,1)$, $F(4)$, $G(3)$, $P(n)$, $Q(n)$, then $L$ is isomorphic to the corresponding Lie superalgebra.}

The hypothesis of this Proposition is obvious if $\gr L=P(n)$ or $A(m,n)$ or $C(n)$: indeed, $L_0=(L_0)_\ev\oplus (L_0)_\od$, where ${(L_0)_\ev= \gr_0 L= L_\ev}$ and $(L_0)_\od =\gr_1 L$. Since $\gr_0 L=L_\ev$ is reductive, any its finite-dimensional representation is completely reducible, in particular, the $L_\ev$-action on~$L_\od$.

The key point in the proof of \textbf{Proposition 2.1.4} is \textbf{Proposition 2.1.3}, which claims that for the Lie superalgebras in question, the  $L_\ev$-module $S^2(L_\od)$ contains $L_\ev$ with multiplicity~ 1.

The cases where $\gr L=\widetilde H(n):=\fh(0|n)$,  or $H(n)=\fh'(0|n)$ are considered in \textbf{Proposition~ 3.3.7, p.~ 66}. Since $\fh(0|n)$ is not simple, there remains only~$H(n)$. 

Finally, the case of $\gr L=S(n)$ is dealt with by \textbf{Theorem 3.3.5, p. 64}. In both these cases, the following two arguments are key ones. 

The first is \textbf{Proposition 3.1.3, p.59} which claims that if 
${L=L_{-1}\supset L_0\supset \dots}$ is a~transitive filtered Lie superalgebra such that $L_0\supset L_\ev$ and ${\dim L_{-1}/L_0=(0|n)}$, then $L$ can be embedded into~$W(n)$. 

The second is the Levi-Malcev theorem for Lie algebras claiming that any semisimple subalgebra of the Lie algebra $L$ can,  by means of conjugations, be embedded into a~fixed Levi subalgebra (i.e., the complement to the radical) of~$L$. Hence, in the case 
where $\gr L=\widetilde H(n)$ or $H(n)$, the Lie superalgebra $L$ contains $\fo(n)$, whereas in the case where $\gr L=S(n)$, the Lie superalgebra $L$  contains~$\fsl(n)$. The rest is a~routine application of the representation theory of simple Lie algebras.
\end{proof}

\smallskip

\textbf{Acknowledgments}. I am thankful to J.~Bernstein, D.~Fuchs, P.~Grozman,   A.~Lebedev, Yu.~Manin, I.~Shchepochkina, and lately B.~Kruglikov and A.~Tikhomirov for education and help; to the grant AD 065 NYUAD for financial support.

\iffalse \AO{suggestion still without my answer: nineteen \texttt{bibitem}s below are never cited ---
\texttt{Car}, \texttt{FL}, \texttt{FiFu}, \texttt{Fr}, \texttt{FrK}, \texttt{Gr1},
\texttt{Gr2}, \texttt{Gr3}, \texttt{K3}, \texttt{K10}, \texttt{Ko}, \texttt{L2}, \texttt{L8},
\texttt{LS1}, \texttt{LSoS2}, \texttt{O}, \texttt{Po1}, \texttt{Ri}, \texttt{Se3}. Journals
usually ask for such a~list to be pruned; you may well prefer to cite some of them instead.
Which is it to be? --- Correcting myself: in my notes of 3 August I~wrote that \texttt{[OSS]},
\texttt{[Pal]} and \texttt{[PalI]} have no \texttt{bibitem} and print as \lq\lq [?]''. That was
wrong: all three are here \textup{(the first and the third are typed
\texttt{\textbackslash bibitem [OSS]}, with a~space, which my check did not match)}, the
bibliography has no dangling entry whatsoever, and the data in them agree with what
I~verified. My apologies for the false alarm.}\fi

\end{document}